\documentclass[11pt]{article}
\usepackage{amsmath,amsfonts, amssymb,graphicx,geometry,authblk,color,soul,comment,hyperref} 
\newcommand{\bfv}{{\mathbf v}}
\newcommand{\bff}{{\mathbf f}}
\newcommand{\bfn}{{\mathbf n}}
\newcommand{\bfzero}{{\mathbf 0}}
\newcommand{\bfK}{{\mathbf K}}
\newcommand{\bfu}{{\mathbf u}}
\newcommand{\bfx}{{\mathbf x}}
\newcommand{\bfe}{{\mathbf e}}
\newcommand{\bfD}{{\mathbf D}}
\newcommand{\bfy}{{\mathbf y}}

\usepackage[dvipsnames]{xcolor}
\usepackage{float}
\usepackage{multirow}

\usepackage{amsthm}

\title{A learning-based multiscale model for reactive flow in porous media}
\author{Mina Karimi and Kaushik Bhattacharya}
\affil{California Institute of Technology, Pasadena CA 91125}

\begin{document}
\maketitle

\begin{abstract}

We study solute-laden flow through permeable geological formations with a focus on advection-dominated transport and volume reactions. As the fluid flows through the permeable medium, it reacts with the medium, thereby changing the morphology and properties of the medium; this in turn, affects the flow conditions and chemistry. These phenomena occur at various lengths and time scales, and makes the problem extremely complex. Multiscale modeling addresses this complexity by dividing the problem into those at individual scales, and systematically passing information from one scale to another. However, accurate implementation of these multiscale methods are still prohibitively expensive. We present a methodology to overcome this challenge that is computationally efficient and quantitatively accurate.  We introduce a surrogate for the solution operator of the lower scale problem in the form of a recurrent neural operator, train it using one-time off-line data generated by repeated solutions of the lower scale problem, and then use this surrogate in application-scale calculations.   The result is the accuracy of concurrent multiscale methods, at a cost comparable to those of classical models.   We study various examples, and show the efficacy of this method in understanding the evolution of the morphology, properties and flow conditions over time in geological formations.
\end{abstract}

\section{Introduction}

The transport of water and solutes through permeable geological formations couple various phenomena ~\cite{lichtner2018reactive, baqer2022review}.  As it flows through the permeable medium, solute-laden water reacts with the medium changing the morphology and properties of the medium; this in turn affects the flow and chemistry.  Permeable medium supports a menagerie of microbial life, and these are affected by and affect the flow and the chemistry.  Chemical interactions, biological activity and flow may cause mechanical failure which in turn changes morphology thereby affecting the flow.  Further, these varied phenomena occur, and manifest themselves, at various length and time scales (see for example \cite{medici_2021} and the citations there): at the {\it molecular} scale where chemistry happens, at the {\it pore} scale where morphology, biology, mechanical properties and flow interact, the {\it core} scale where averaged features emerge, and finally at the {\it geological} scale relevant to aquifer and reservoir hydrology.  Furthermore, these geological formations are heterogenous at various scales ranging from the individual grains, to pores, to formations, fissures, caverns and cracks to geological strata and faults.   Finally, they are highly anisotropic at various scales ranging from the clay particles to faults to geological formations.

This enormous complexity makes modeling extremely difficult.   There are well-developed models for individual processes at particular scales, and in some cases conceptual frameworks to link various scales.  Still, {\it there is a critical gap} in accurate methods of passing information from one scale to another, especially when it concerns multiple phenomena and history (or time)-dependent phenomena.  In this paper, we address this coupling between the pore, core and geological formation scales with a focus on advection-dominated transport and volume reactions.   

We begin with a brief review of the existing approaches, and then describe our contributions.

\noindent {\it Geological scale models.} Continuum models for reactive flow at the geological scale focus on predicting the overall flow and transport. These models usually employ empirical relations to relate the evolving transport parameters to changes in porosity, reactivity, and specific area \cite{baqer2022review}. To do so, many studies estimate the change in porosity \cite{appelo2004geochemistry}, diffusivity \cite{seigneur2019reactive}, permeability \cite{appelo2004geochemistry, witherspoon1980validity, oron1998flow}, reaction rate, and specific area \cite{steefel1998multicomponent} due to chemical reactions. However, these approaches face some limitations.  These models are derived assuming simple geometries (e.g. spherical grains) and therefore can not provide detailed insight into the pore structure and specific area, and their evolution.  Further, they have a limited range of accuracy in the presence of precipitation and dissolution.  In particular, they often fail when porosity approaches critical values or percolation.

\noindent {\it Core scale models.} Core-scale simulations focus on predicting the evolution of the pore structure and average fluid transport properties. Using conventional finite element methods for these simulations is computationally expensive due to the need to update the geometry \cite{baqer2022review}. Several alternative numerical strategies have been developed to overcome this limitation, such as pixel-based approaches \cite{kang2003simulation, kang2002lattice}, level set methods \cite{li2008level, varloteaux2013reactive}, smoothed particle hydrodynamics \cite{tartakovsky2016smoothed}, and adaptive discontinuous Galerkin finite element methods \cite{sun2006dynamic}.

\noindent {\it Multiscale models.} These models employ upscaling approaches to bridge the gap between pore and geological scales. Volume averaging methods \cite{bear2012introduction} estimate the effective transport properties and derive macroscopic equations. Initially explored for porous media dispersion by \cite{carbonell1983dispersion, quintard1993transport}, these methods were later extended to spatially periodic porous media for reactive fluid transport by Paine \cite{paine1983dispersion}. However, such an approach necessitates incorporating ad hoc assumptions for closure relations.

Early studies by Auriault and Adler \cite{auriault1995taylor} and Rubinstein and Mauri \cite{rubinstein1986dispersion} investigate multiscale expansions to upscale Taylor dispersion \cite{taylor1953dispersion} in periodic porous media, primarily focusing on diffusion-dominated, non-reactive transport. Later, this methodology was extended to reactive flow with surface reactions by Mauri \cite{mauri1991dispersion}. However, this approach failed to derive effective properties for advection-dominated flows. Other researchers \cite{maruvsic2005homogenization, allaire2007homogenization, donato2005averaging} introduced a two-scale expansion approach with drift to address advection-dominated flows. They employed a moving coordinate system with an effective drift velocity to investigate advection's impact on effective equations. This method was further applied by Allaire et al., \cite{allaire2010two} to study reactive transport with surface reactions while considering a linear adsorption/desorption rate at the solid interface.

\noindent {\it Machine learning.} In recent years, machine learning (ML) techniques have been applied to expedite reactive transport simulations and estimate effective porous medium properties. Liu et al \cite{liu2022machine} have used an ML approach to predict the effective reaction rate using features of pore structures. Artificial neural networks have been used to accelerate the geological scale reactive flow simulations \cite{demirer2022improving, sprocati2021integrating}, and deep learning methods have been used to estimate the permeability considering non-reactive transport in the porous medium \cite{wang2022pore, wang2021physics}. These consider a snap-shot in time and do not address evolution over long periods. Wang and Battiato \cite{wang2021upscaling} have developed a deep-learning multiscale model to predict the clogging caused by solute precipitation in a microcrack network. Lu et al. \cite{lu2021data} have used a neural network model to predict the evolution of the uranium distribution coefficient in the subsurface due to thermal and chemical processes. 

In the context of mechanical properties, recurrent neural networks (RNNs) address the history-dependent behaviors \cite{medsker2001recurrent}, Long Short-Term Memory (LSTM) effectively remembers information over long time intervals \cite{ghavamian2019accelerating}, gated recurrent unit (GRU), a simplified variant of LSTM captures similar temporal relationships in data \cite{mozaffar2019deep, wu2020recurrent}.   However, LSTM-based approaches require millions of variables to be trained from data.  Recently, a recurrent neural operator (RNO) has been introduced inspired by internal variable theory, offering an efficient approach for multiscale modeling \cite{liu2023learning, liu2022learning}.

\paragraph{Our Contributions.} In this work, we develop a methodology for investigating flow through underground geological formations over long periods of time. We specifically focus on an advection dominated transport in a porous medium and reactions at the fluid-solid interface characterized by nonlinear reaction kinetics.  We first use a two-scale expansion method with effective drift velocity \cite{allaire2010two} to obtain the governing equations are the core and geological formulations.  This is described in Section \ref{sec:2scalemodel}.   

We then introduce a recurrent neural operator (RNO), and use it to learn the solution operator of the core scale problem.  To elaborate, we use data generated by repeated solutions of the core scale problem to train the RNO to learn the map between geological scale variables over time (e.g., velocity and solute concentration histories) and the effective transport properties (e.g., permeability, diffusivity, advection velocity, porosity, and specific area). This is described in Section \ref{sec:RNO}.

Finally, we use the trained RNO surrogate model as a surrogate in geological scale computations to investigate the long time evolution of these formations in the presence of volume reactions induced by flow and transport.  We demonstrate the accuracy of the approach, as well as its ability to reveal non-trivial interactions between the core and geological scale in Section \ref{sec:Multiscale simulation}.

We conclude in Section \ref{sec:conc} with a discussion of implications, promises and open issues.

\section{Two scale model} \label{sec:2scalemodel}

\subsection{Governing equation and non-dimensionalization} \label{sec:gov}
We consider a porous geological formation $\tilde \Omega$ composed of a porous region $\tilde \Omega_p$ and solid region $\tilde \Omega_s$ with $\tilde \Omega_p \cap \tilde \Omega_s = \phi, \ \tilde \Omega^c = \tilde \Omega_p^c\cup \tilde \Omega_s^c$ where the superscript `c' denotes the closure.  We have an incompressible fluid flow in the pore governed by the steady Stokes equation
\begin{align} \label{eq:stokes}
    \begin{cases}
        -\tilde{\nabla}\tilde{p} + \tilde{\nu}\Delta \tilde{\bfv} + \tilde{\bff} = \bfzero \quad & \text{in} \ \tilde{\Omega}_p, \\
        \tilde{\nabla}\cdot\tilde{\bfv} = 0 & \text{in} \  \tilde{\Omega}_p, \\
        \tilde{\bfv} = \bfzero & \text{in} \ \tilde{\Omega}_s
    \end{cases}
\end{align}
where $\tilde{\bfv}$ is the particle velocity, $\tilde{p}$ is pressure, $\tilde{\nu}$ is the viscosity and $\tilde{\bff}$ is the body force.  We assume that the fluid carries a solute that is transported in the fluid through a combination of diffusion and advection and reacts with the surface of the solids, resulting either in deposition or dissolution
\begin{align} \label{eq:transport}
\begin{cases}
    \tilde{c}_{\tilde{t}} + \tilde{\bfv}\cdot\tilde{\nabla}\tilde{c} -  D \tilde{\Delta}\tilde{c} = 0   \quad & \text{in} \ \tilde{\Omega}_p, \\
    -D\tilde{\nabla}\tilde{c}\cdot\bfn = -(\tilde{c}- \tilde{m})\tilde{v}_n = \tilde{q}_c(\tilde c) & \text{on} \ \partial\tilde{\Omega}_p
\end{cases}
 \end{align}
where $\tilde c$ is the concentration of the solute, $D$ is the diffusion constant, $\bfn$ the unit normal vector to the surface $\partial\tilde{\Omega}_p$ oriented outward with respect to $\tilde{\Omega}_p$, $\tilde m$ is the concentration of the solute in the solids, $\tilde{v}_n$ is the normal speed of the solid/fluid interface due to dissolution or deposition and $\tilde{q}_c(c)$ is the reaction rate that depends on the solute concentration.  Above, the second term of the bulk equation describes advection, while the third describes diffusion.  On the interface, the first equation describes the mass balance between the flux of the solute from the fluid to the interface and the growth of the interface, while the second relates the growth of the interface to the interfacial reaction rate. 
These equations have to be supplemented with appropriate boundary conditions.

Note that we have assumed that the fluid flow is steady while the solute transport is time-dependent.  We assume that the reaction rate at the interface and, consequently, the rate of reconstruction of the porous medium is slow compared to the time-scale of the fluid flow.  So, the fluid flow reaches a steady state at each time as the medium reconstructs.

Now, the characteristic length of the geological formation $L$ is very large compared to the characteristic length of the core or representative volume $\ell$, $L >> \ell$.  This makes the system (\ref{eq:stokes}, \ref{eq:transport}) difficult to solve: we have to resolve the flow and transport with a resolution small compared to $\ell$, but on a domain of size $L$.  Therefore, we resort to a two-scale asymptotic expansion under the assumption that the ratio of length-scales $\epsilon = \ell/L$ is small, $\epsilon<<1$.

In order to do so, we change to non-dimensional units by scaling length with the characteristic length $L$ of the geological formation, the velocity with the characteristic velocity $V$ and the pressure with characteristic pressure $\Pi$.  It follows that time is rescaled by the characteristic time $T=L/V$.  We expect slow flows with small $V$ over long distances $L$, which means that the characteristic time $T$ is large and consistent with steady state.  Now, recall that the characteristic length of the pores is small, and therefore, in order to have non-trivial flow, we need the viscosity to be extremely small.  Therefore, we assume that the characteristic viscosity is $\Pi T/\epsilon^2$.  The non-dimensional flow equations are given by
\begin{align}\label{eq: dimensionless stokes}
    \begin{cases}
        -{\nabla}{p^\epsilon} + \epsilon^2 \nu \Delta {\bfv^\epsilon} + {\bff} = \bfzero \quad & \text{in} \ {\Omega}^\epsilon_p, \\
        {\nabla}\cdot{\bfv^\epsilon} = 0 & \text{in} \  {\Omega}^\epsilon_p, \\
        {\bfv^\epsilon} = \bfzero &\text{in} \ {\Omega}^\epsilon_s.
    \end{cases}
\end{align}
Above, we have introduced $\epsilon$ as a superscript in the non-dimensional variables to signify that the porosity and therefore variations in these quantities are at a scale $\epsilon$.

We now turn to solute transport.  We non-dimensionalize the concentration with a characteristic concentration $C$ and introduce two non-dimensional numbers, the P\'eclet number Pe and the Damk\"ohler number Da:
\begin{equation} \label{eq:pecdam}
    \text{Pe}^\epsilon = \frac{LV}{D}, \quad \text{Da}^\epsilon = \frac{LK}{D} 
\end{equation}
where $K$ is a characteristic reaction rate.   The non-dimensional equations of solute transport are
\begin{equation}\label{eq: dimensionless concentration}
    \begin{cases}
        c_t^\epsilon + \text{Pe}^\epsilon\bfv\cdot \nabla c^\epsilon - \Delta c^\epsilon = 0 \quad & \text{in} \ \Omega^\epsilon_p \\
      -\nabla c^\epsilon \cdot \bfn = - v_n^\epsilon (c^\epsilon-m) = \text{Da}^\epsilon q_c^\epsilon & \text{on} \ \partial \Omega^\epsilon_p
    \end{cases}
\end{equation}

Now, we are interested in situations where we have significant advection at the pore scale, and very slow reactions at the interface.  So, we assume that
\begin{equation} \label{eq:pecdamscale}
    \text{Pe}^\epsilon = \frac{\widehat{\text{Pe}}}{\epsilon}, \quad \text{Da}^\epsilon =  \epsilon \widehat{\text{Da}}, \quad v_n^\epsilon = \epsilon \hat{v}_n
\end{equation}
where $\widehat{\text{Pe}}, \widehat{\text{Da}}, \hat{v}_n$ are all $O(1)$.

\subsection{Two-scale model} \label{sec:2scalesum}

We assume that the porous media is almost periodic, i.e., it is periodic on the scale $\epsilon$ but can change over long distances compared to $\epsilon$.  To be precise, we assume $\Omega^\epsilon_p (x) = \Omega_p(x,x/\epsilon)$ where $\Omega_p$ is periodic with period 1 in the second variable; so porous medium is periodic on $\epsilon Y$ where $Y$ unit cube or unit cell in the vicinity of the point $\bfx$.  Further, $Y = Y_p (\bfx) \cup Y_s (\bfx)$ where $Y_p(\bfx)$ is the pore in the unit cell in the vicinity of the point ${\bf x}$ in the geological formation. We show that under this assumption, we can approximate the solution of the system (\ref{eq: dimensionless stokes},\ref{eq: dimensionless concentration}) by solving a geological scale problem where the constitutive behavior is determined by solving a core scale problem.  We first describe the two problems, and then justify this derivation in the following subsections.

\subsubsection{Geological scale model}
We can find the overall flow and solute transport at the geological scale by solving the following system on $\Omega$:
\begin{align} \label{eq:reservoir}
\begin{cases}
\nabla\cdot \bfv_0 = 0, \quad \bfv_0 = \frac{1}{\nu} \bfK^* (\bff - \nabla p_0),\\
\lambda {c}_{0t} +  \lambda \text{Pe} \overline \bfv \cdot \nabla c_0 - \nabla\cdot (\bfD^* \nabla c_{0})  = - \gamma \widehat{\text{Da}} q_c
\end{cases}
\end{align}
for the overall velocity $\bfv_0 $, pressure $p_0$, concentration $c_0$ at the geological scale subject to boundary conditions.  Above, the parameters $\bfK^* = \bfK^*(t,\bfx)$ is the permeability tensor,  $\bfD^* = \bfD^*(t,\bfx)$ is the effective diffusivity tensor, $\lambda = \lambda(t,\bfx)$ is the pore volume fraction, $\overline \bfv=\overline \bfv(t,\bfx)$ is an effective advection velocity, and $\gamma = \gamma(t,\bfx)$ is a local surface area per unit volume.  Note that these parameters are all functions of time, and they are specified through the core scale problem below.

Note that this geological scale problem is solved at the entire geological $\Omega$ and the coefficients vary only on the scale of the geological formation.  All the information about the pores has been subsumed into parameters that only vary at the geological scale.

\subsubsection{Core scale or unit cell model} \label{sec:core}
Given a porous unit cell at the macroscopic point $\bfx$ at time $t$, i.e., given $Y_p$, the unit cell problem is to solve
\begin{align} \label{eq: unit cell}
    \begin{cases}
            -\nabla_y q^i + \nu \Delta_y \bfu^i + \bfe^i = \bfzero, \ \  \nabla_y \cdot \bfu^i = 0 \quad & \text{in} \ Y_p \\
	    \bfu^i = \bfzero & \text{in} \ Y_s\\
	    \widehat{\text{Pe}} \ \bfv\cdot \left( \nabla_y \chi^j \right)  - \Delta_y \chi^j  = \
	    (\bfv^* - \widehat{\text{Pe}} \ \bfv)\cdot\bfe^j \quad & \text{in} \ Y_p\\
       - \left( \nabla_y \chi^j \right)\cdot {\bfn} = \bfe^j\cdot {\bfn} & \text{on} \ \partial Y_p
    \end{cases}
\end{align}
periodic velocity fluctuation $\bfu^i$, pressure fluctuation $q^i$, chemical fluctuation $\chi^j$ when the overall flow is in the direction $\bfe^i$ and overall solute transport in the direction $\bfe^j$ for $i,j = 1, \dots, d$ (dimension $d$).  Above, we use $\nabla_y, \Delta_y$ to signify that these are derivatives with respect to the spatial variable in the unit cell.

We can then use it to find the parameters 
\begin{align}
\bfK^*_{ij} &= \int_{Y_p} \nabla_y \bfu^i \cdot \nabla_y \bfu^j \ dy, \label{eq:K}\\
\bfD^*_{ij} &= \int_{Y_p} \bfe^i\cdot \bfe^j \ dy + \widehat{\text{Pe}} \int_{Y_p} (\overline \bfv_i -  \bfv_i) \chi^j \ dy 
    +  \int_{Y_p} \nabla_y \chi^j\cdot \bfe^i \ dy , \label{eq:D}\\
\overline \bfv &= \frac{1}{|Y_p|} \int_{Y_p} \bfv \ dy, \quad 
\lambda = |Y_p|, \quad 
\gamma = |\partial Y_p|. \label{eq:v}
\end{align}

Finally, the microstructure evolves according to the equation 
\begin{equation}\label{eqn: interface speed}
    \hat{v}_n = \frac{\widehat{\text{Da}}}{ m - \tilde c_0} q_c(\tilde c_0) \approx \frac{\widehat{\text{Da}}}{m}  q_c(\tilde c_0).
\end{equation}
The last approximation uses the fact that $\tilde c_0 << m$.  This evolution happens on the geological time scale, and therefore the parameters above change on the geological time scale.

\subsubsection{Summary}

\begin{figure}
\centering
    \includegraphics[width=0.6\textwidth]{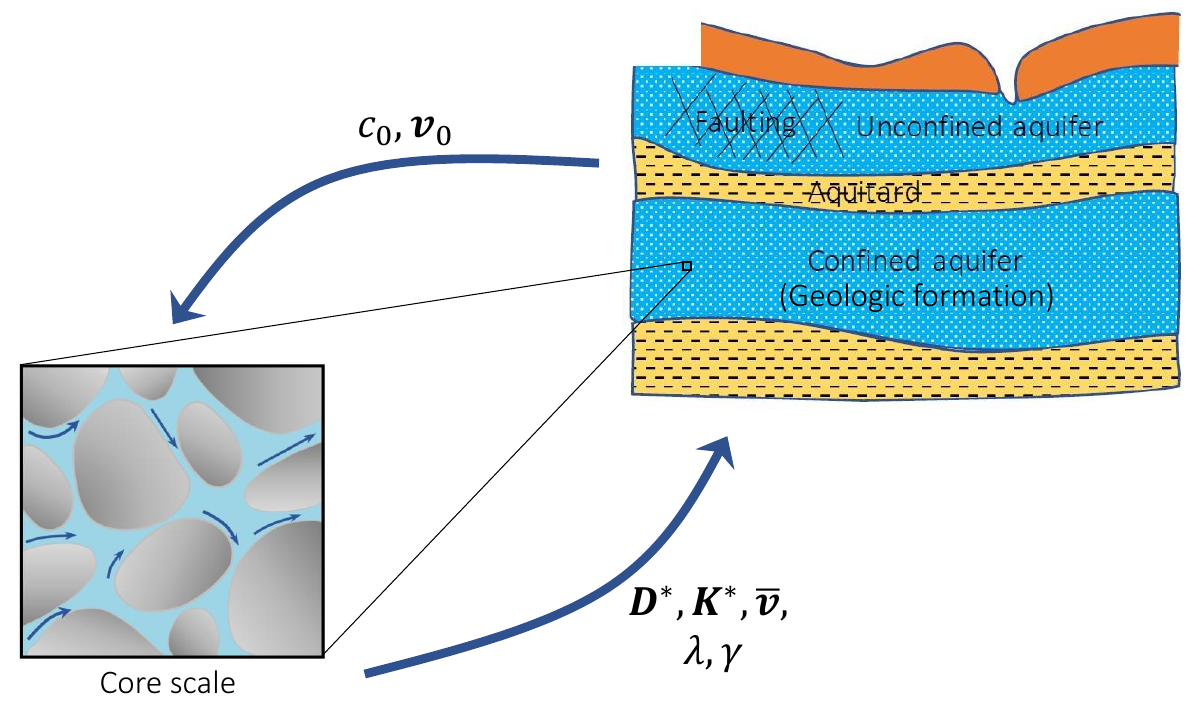}
    \caption{Schematic figure of the two-scale framework.}
    \label{fig:Sec2_Fig1}
\end{figure}

We summarize the resulting multiscale formulation.  We solve (\ref{eq:reservoir}) on the geological scale.  This requires us to obtain the parameters $\bfK^*(t,\bfx), \bfD^*(t,\bfx), \bar \bfv (t,\bfx), \gamma(t,\bfx), \lambda(t,\bfx)$.  To do so, at each point $\bfx$ at the geological scale, we provide a history of the flow $\bfv_0$ and solute concentration $c_0$ to a unit cell, solve (\ref{eq: unit cell}) and obtain the parameters from (\ref{eq:K}), (\ref{eq:D}) and (\ref{eq:v}) as a function of time while evolving the microstructure according to (\ref{eqn: interface speed}).  This is illustrated in Figure \ref{fig:Sec2_Fig1}.

While this formulation separates the original problem into two separate problems, it is still computationally demanding: we have to solve a core scale problem at at every quadrature point of the geological formation and at every instant of time.  The direct implementation of this framework is often referred to as concurrent multiscale approach.  We propose an alternate approach in Sections \ref{sec:RNO} and \ref{sec:Multiscale simulation} below.

\subsection{Details of the asymptotic expansion}

This sub-section outlines the derivation of the two-scale formulation above.  It closely follows Allaire et al.~\cite{allaire2010two}.

\subsubsection{Porous media flow}

We look for a solution to the system (\ref{eq: dimensionless stokes}) with the ansatz
\begin{align}\label{eq: stokes ansatz}
    \bfv^{\epsilon} (\bfx) = \sum_{i=0}^{+\infty} \epsilon^i \bfv_i \left(\bfx, \frac{\bfx}{\epsilon} \right) \, \quad p^{\epsilon} (\bfx) = \sum_{i=0}^{+\infty} \epsilon^i p_i \left( \bfx, \frac{\bfx}{\epsilon} \right).
\end{align}
By collecting terms and then using the Fredholm alternative, L\'evy \cite{levy1983fluid} showed that $\bfv_0 = \bfv_0 (\bfx), \ p_0 = p(\bfx)$ (i.e, they are independent of the fast variable), and that these satisfy (\ref{eq:reservoir})$_1$ with $\bfK^*$ given by (\ref{eq:K}) where $\bfu^i$ satisfies (\ref{eq: unit cell})$_{1,2}$.  Further, up to leading order the velocity and the pressure are
\begin{align} \label{eq:velpre}
 \bfv^{\epsilon} (\bfx) = \bfv_0(\bfx) + \epsilon \sum_{i=1}^d (\bfv_0)_i(\bfx) \bfu^i \left(\frac{\bfx}{\epsilon}\right), \quad 
 p^{\epsilon} (\bfx) = p_0(\bfx) + \epsilon  \sum_{i=1}^d q^i \left(\frac{\bfx}{\epsilon}\right) \left( \bff_i(\bfx) - \frac{\partial p_0}{\partial x_i}(\bfx) \right) .
 \end{align}

\subsubsection{Solute transport}
We look for a solution to the system (\ref{eq: dimensionless concentration}) with the ansatz
\begin{align}
    c^\epsilon (t, \bfx) = \sum^{+\infty}_{i=0} \epsilon^i \tilde c_i \left( t, \bfx-\frac{\bfv^*}{\epsilon}t , \frac{\bfx}{\epsilon} \right)
\end{align}
where $ \bfv^*$ is a {\it drift velocity} to be determined.   We have chosen a macroscopic coordinate $\bfx' = \bfx-\frac{\bfv^*}{\epsilon}t$ that is not stationary in the geological formation in order to account for the advection.  We substitute the ansatz into (\ref{eq: dimensionless concentration}), apply the chain rule, 
\begin{align}
    c_t = \frac{\partial{\tilde c}}{\partial t} - \frac{\bfv^*}{\epsilon}\cdot \nabla' \tilde c , \quad
    \nabla c = \nabla' \tilde c + \frac{1}{\epsilon}\nabla_y \tilde c
\end{align}
and collect terms with different powers of $\epsilon$.

The leading order (smallest power) is $\epsilon^{-2}$, and we have the following equations at that order:
\begin{equation}
    \begin{cases}
        \hat{\text{Pe}} \bfv\cdot\nabla_y \tilde c_0 - \Delta_y \tilde c_0 = 0 \quad &  \text{in} \ Y_p,\\
      -\nabla_y \tilde c_0 \cdot {\bfn} = 0 &\text{on} \ \partial Y_p
    \end{cases}
\end{equation}
which yields to $\tilde c_0= \tilde c_0(\bfx')$.

At the next order, $\epsilon^{-1}$, we have 
\begin{equation}\label{eq: concentration order -1}
    \begin{cases}
        \widehat{\text{Pe}} \ \bfv\cdot  \nabla_y \tilde c_1  - \Delta_y  \tilde c_1  =  (\bfv^* - \widehat{\text{Pe}} \ \bfv ) \cdot\nabla' \tilde c_0 \quad & \text{in} \ Y_p,\\
      -\left(  \nabla_y \tilde c_1 \right)\cdot {\bfn} = \nabla' \tilde c_0\cdot {\bfn}  & \text{on} \ \partial Y_p.
    \end{cases}
\end{equation}
We seek to solve this equation for $\tilde c_1$.  This is possible if and only if the following compatibility condition (Fredholm alternative) is satisfied.
\begin{equation}
    \int_{Y_p}  (\bfv^* - \widehat{\text{Pe}})  \cdot\nabla' \tilde c_0 \ dy - \int_{\partial Y_p} \left( \nabla' \tilde c_0\cdot {\bfn} \right) d S = 0 
    \quad  \Rightarrow \quad
    \bfv^* = \frac{\widehat{\text{Pe}}}{|Y_p|} \int_{Y_p} \bfv \ dy.
\end{equation}
Recalling the expansion  (\ref{eq:velpre})$_1$ for $\bfv$, we can re-write the drift velocity as follow
\begin{equation}
    \bfv_i^* = \frac{\widehat{\text{Pe}}}{|Y_p|} \sum_{i=1}^d (\bfv_0)_i \int_{Y_p} \bfu^i \ dy
\end{equation}
Returning to (\ref{eq: concentration order -1}), we notice that the solution $\tilde c_1$ will depend linearly on $\nabla' c_0$.  So, set
\begin{equation} \label{eq: c1}
    \tilde c_1(t, \bfx', \bfy) = \sum^d_{i=1} \frac{\partial \tilde c_0}{\partial x'_i} (t, \bfx') \chi^i(\bfy).
\end{equation}
It follows that $\chi^i$ satisfy (\ref{eq: unit cell})$_{3,4}$.

Now, turning to order $\epsilon$, we have
\begin{equation}
    \begin{cases}
        \widehat{\text{Pe}} \ \bfv \cdot \nabla_y \tilde c_2 - \Delta_y c_2 \\
        \quad \quad = - \frac{\partial \tilde c_0}{\partial t} + (\bfv^*  - \widehat{\text{Pe}} \ \bfv )\cdot \nabla' \tilde c_1 
        + \nabla' \cdot\left( \nabla' \tilde c_0 + \nabla_y \tilde c_1 \right) + \nabla_y \cdot(\nabla' \tilde c_1)  \quad & \text{in} \ Y_p,\\
      -\nabla_y \tilde c_2\cdot{\bfn} - \nabla' \tilde c_1\cdot\bfn =  \hat{v}_n (\tilde c_0 - m) = \widehat{\text{Da}} \ q_c ( \tilde c_0) & \text{on} \ \partial Y_p.
    \end{cases}
\end{equation}
This equation has a solution for $\tilde c_2$ if and only if the following compatibility condition (Fredholm alternative) is satisfied
\begin{equation}
    \int_{Y_p} \left( - \frac{\partial \tilde c_0}{\partial t} + (\bfv^*  - \widehat{\text{Pe}} \ \bfv )\cdot \nabla' \tilde c_1 
        + \nabla' \cdot\left( \nabla' \tilde c_0 + \nabla_y \tilde c_1 \right)   \right) \ d y = \int_{\partial {Y_p}}\widehat{\text{Da}} \ q_c(\tilde c_0 ) \ d S .
\end{equation}
Substituting for $\tilde c_1$, we obtain the homogenized equation
\begin{equation}\label{eqn: homogenized concentration equation}
    -\lambda \frac{\partial \tilde c_0}{\partial t}  = -\nabla' \cdot (\bfD^* \nabla' \tilde c_0) + \gamma \widehat{\text{Da}} \ q_c (\tilde c_0)
\end{equation}
where $\lambda = |Y_p|$, and $\gamma = |\partial Y_p|$, and the effective diffusion tensor $\bfD^*$ is given by (\ref{eq:D}).
We then transform back to stationary coordinates by setting.  
\begin{equation}
c_0(t, \bfx) = \tilde c_0 \left( t, \bfx-\frac{\bfv^*}{\epsilon}t  \right)
\end{equation}
to obtain (\ref{eq:reservoir})$_2$.

Finally, the boundary condition at $O(\varepsilon)$ gives (\ref{eqn: interface speed}).

\section{Learning the core scale behavior} \label{sec:RNO}

The two-scale formulation above requires us to solve the core scale problem at each time step and each quadrature point in the geological formation.  This is prohibitively expensive.  So, we seek to ``learn'' the solution operator of the core scale problem.  Specifically, for a given initial microstructure $T^0$, we view the core scale problem as a map from the velocity and concentration history to the current permeability, diffusivity, drift velocity, specific area, and pore volume fraction.  
\begin{equation}
\Phi: {\mathcal I}[0,t] \to {\mathcal O}(t), \quad {\mathcal I}(\tau) = \{\bfv_0({\tau}), c_0 (\tau)\}, \ {\mathcal O}(\tau) = \{ \bfK^*(\tau), \bfD^*(\tau), \bar \bfv (\tau), \gamma(\tau), \lambda(\tau)\}.
\end{equation}
where the input ${\mathcal I}[0,t]= \{\bfv_0(\tau), c_0 (\tau): \tau \in [0,t] \}$ is specified over the time interval $[0,t]$ and the output ${\mathcal O}(t) = \{ \bfK^*(t), \bfD^*(t), \bar \bfv (t), \gamma(t),  \lambda(t)\}$ at time $t$. 
We seek an  approximation in the form of a parametrized map
\begin{equation}
\Psi: {\mathcal I}[0,t] \times {\mathbb R}^p \to {\mathcal O} 
\end{equation}
and train it using data $\{ {\mathcal I}^n, {\mathcal O}^n \}_{n=1}^N$ that is generated using numerical simulation of $\Phi$.  In other words, we shall postulate a form for $\Psi$ and find the parameters $\Theta^*$ that minimizes a loss function for data generated by repeated solutions of the core-scale problem.

\subsection{Recurrent neural operator}
There are two issues we have to address in postulating an approximation $\Psi$.  First, the map $\Phi$ (and hence $\Psi$) has as its input a function defined on an interval of time.  Thus our approximation has to be an operator.  One idea is to discretize the functions in time and then find a neural network approximation with the discretized function as input.   Unfortunately, this approximation will depend on the discretization (time step), and hence can only be used at that discretization.  However, it is natural in a multi-scale setting to use different discretization for the core scale problem (generating data) and the geological scale problem.  Further, one may use an adaptive discretization in the geological scale problem.  For these and other reasons, we want approximation to be independent of discretization.   Second, the output at time $t$ depends on the history of the input.  So, we want our map to be history dependent. 
 
Following experience and practice in continuum physics, we postulate that the history can be encapsulated in $k$ state or internal variables $\{\xi_\alpha\}_{\alpha=1}^k$ that evolve in time.  Then, following recent work in the multi-scale modeling of mechanical properties of materials \cite{bhattacharya_2023,liu2023learning},  we postulate $\Psi$ to be a recurrent neural operator :
\begin{equation} \label{eq:rno}
\Psi:     \begin{cases}
        {\mathcal O}(t) =f\left( {\mathcal I}(t), \{\xi_\alpha(t)\}_{\alpha=1}^k; \Theta \right), \\
        \dot{\xi}_i(t) =g_i \left( {\mathcal I}(t), \{\xi_\alpha(t)\}_{\alpha = 1}^k; \Theta \right), \quad i = 1, \dots, k
    \end{cases}
\end{equation}
where $f, g_i$ are feed-forward deep neural networks parametrized by $\Theta$ (weights and biases).  
The architecture (\ref{eq:rno}) is formulated to be in continuous time.  To implement it with time discretization, we use a backward Euler discretization: 
\begin{equation}
    \begin{cases}
        {\mathcal O}^n = f\left( {\mathcal I}^n, \{\xi^n_\alpha\}_{\alpha=1}^k \right) \\
        \xi^n_i = \xi^{n-1}_i + (\Delta t_n)g_i \left( {\mathcal I}^n, \{\xi^{n-1}_\alpha\}_{\alpha = 1}^k \right), \quad i= 1, \dots, k.
    \end{cases}
\end{equation}
Note that $f, g_i$ and the internal variables, and therefore the approximation, is independent of the discretization.

The number $k$ of internal variables has to be chosen \textit{a priori}, but the actual internal variables are identified from the data as a part of the learning.  As noted, they encapsulate the history dependence.  They do not necessarily have any intrinsic physical meaning.  Indeed, note that the form of the architecture (\ref{eq:rno}) is invariant under the reparametrization $\xi' = \Xi (\xi)$ for any diffeomorphism $\Xi$.   In some special examples, it is possible to choose a parametrization so that the internal variables are interpretable~\cite{liu2023learning}; however, this is not always the case.  We refer the reader to ~\cite{liu2023learning} for a discussion of these and other aspects of this architecture.

\subsection{Data and training}
We generate the data by solving the core scale problem (\ref{eq: unit cell}) over some interval $[0,T]$ to yield our data in the form $\{{\mathcal I}^n[0,t],{\mathcal O}^n[0,T]\}_{n=1}^N$.  To do so, we have to sample the inputs $\{{\mathcal I}^n[0,t]\}_{n=1}^N$ in such a manner that is rich enough to represent the actual trajectories encountered in the geological scale model.  Broadly, we anticipate trajectories of velocity and concentration that can vary over time, and also change slope as some region gets clogged or fully dissolved.  So we use the following strategy.   We take our interval $[0,T]$ and divide them into $M$ sub-intervals $\{[t^{m-1},t^m] \}_{m=1}^M$ with $t^m \le t^{m+1}, \ t^0 = 0, \ t^M=T$ where $\{t^m\}_{m=1}^{M-1}$ are chosen from a uniform distribution (and relabelled to be increasing).  We then set each component ${\mathcal I}_i$ of the input at times $\{t^m\}_{m=1}^{M-1}$
\begin{equation}
   ({\mathcal I}_i)(t^m) = (\mathcal{I}_i)^{m-1} + \nu^m {\mathcal I}_i^{\max} \sqrt{ t }, ~~ \text{with} ~~~ i=1 ~~ \text{for} ~~ c_0, ~~ i=1, \cdots, d ~~ \text{for} ~~ \bfv_0
\end{equation}
where $\nu_n \in \{ -1, 1\}$.  We then obtain ${\mathcal I}_i [0,T]$ via a cubic spline interpolation. We refer the reader to  ~\cite{liu2023learning, liu2022learning} for a discussion.  We clarify that $\{t^m\}$ is distinct from the time steps used for generating the data. We consider $T=1$,  $M=5$ and use 200 time-steps to calculate the solution. 

We emphasize that the data provided to the RNO is at the geological scale: inputs (velocity and concentration history at a point), and output (permeability, diffusivity, advection velocity, specific area, and pore volume fraction). There is no information about the pore scale.  The internal variables are inferred from this data as a part of the training process.

After generating the data, we proceed to train the RNO. The training process involves sequentially feeding the input data into the network, computing outputs at each time step, and comparing them to the target values during the forward propagation. We then use the backpropagation through time
to calculate gradients over the entire sequence, and optimize the parameters.

\subsection{Results}\label{sec:RNO results}

\begin{figure}
\centering
    \includegraphics[width=0.3\textwidth]{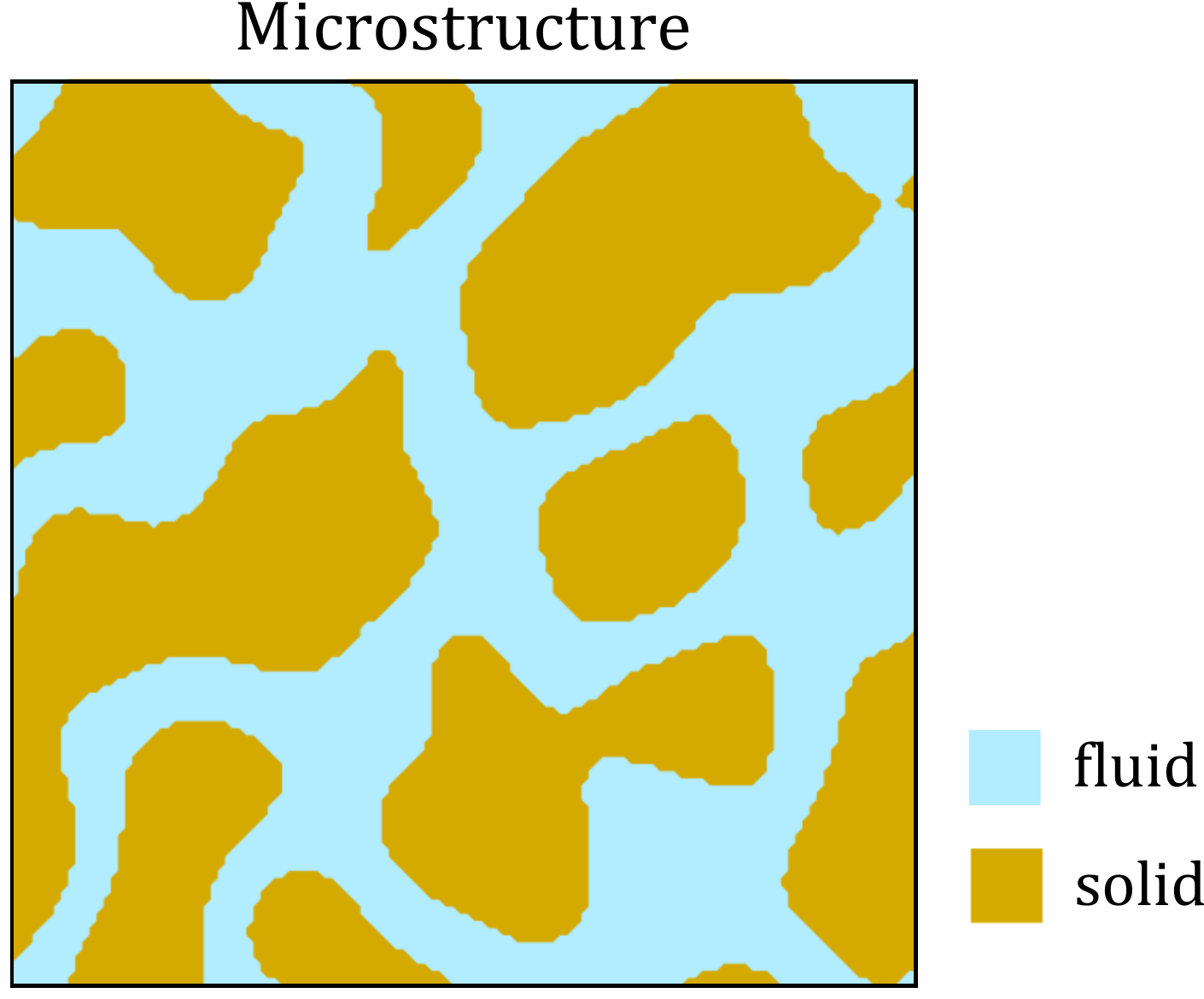}
    \caption{Initial core-scale microstructure}
    \label{fig:micro}
\end{figure}

We now demonstrate the ability of the RNO to approximate the core-scale problem.  We consider an initial microstructure shown in Figure \ref{fig:micro}. We generate our data by solving the core-scale problem described in Section \ref{sec:core} using a numerical algorithm described elsewhere \cite{karimi_2023}.

We consider a fully connected $4$-layer neural network, with each layer consisting of $200$ nodes. We use the nonlinear activation function scaled exponential linear unit (SELU)~\cite{klambauer2017self}, and optimize the parameters using the ADAM optimization algorithm~\cite{kinga2015method}.  We consider the following loss function
\begin{equation} \label{eq:loss}
{\mathcal L} = \frac{1}{D_\text{train}} \sum_{d=1}^{D_\text{train}} 
\frac{\int^T_0\left|  | \overline{\mathcal{O}}^{\text{truth}}_d -  \overline{\mathcal{O}}^{\text{approx}}_d \right|^2 d t}{\int^T_0\left|  | \overline{\mathcal{O}}^{\text{truth}}_d \right|^2 d t}, 
\end{equation}
where $d$ indexes the trajectory in the training dataset and $\overline{\mathcal{O}}$ is the normalized output.  To compute this, we use min-max normalization on  each physical component of the output, 
\begin{equation}
(\overline{\bfD}^*_{ij})_d = \frac{ (\bfD^*_{ij})_d - \min{p,q,r} (\bfD^*_{pq})_r}{\max{p,q,r} (\bfD^*_{pq})_r-\min{p,q,r} (\bfD^*_{pq})_r}
\end{equation}
and so forth.  We define and train the RNO with $\log \bfK^*$ instead of $\bfK^*$ to properly emphasize the almost clogged regime.

Recall that we have to fix the number of internal variables $k$ before training.  We do so, repeating the training for $k = 0, 1, \dots, 4$.  Thus we have five fully trained RNOs with differing numbers of internal variables.

\begin{figure}[t]
\centering
    \includegraphics[width=4in]{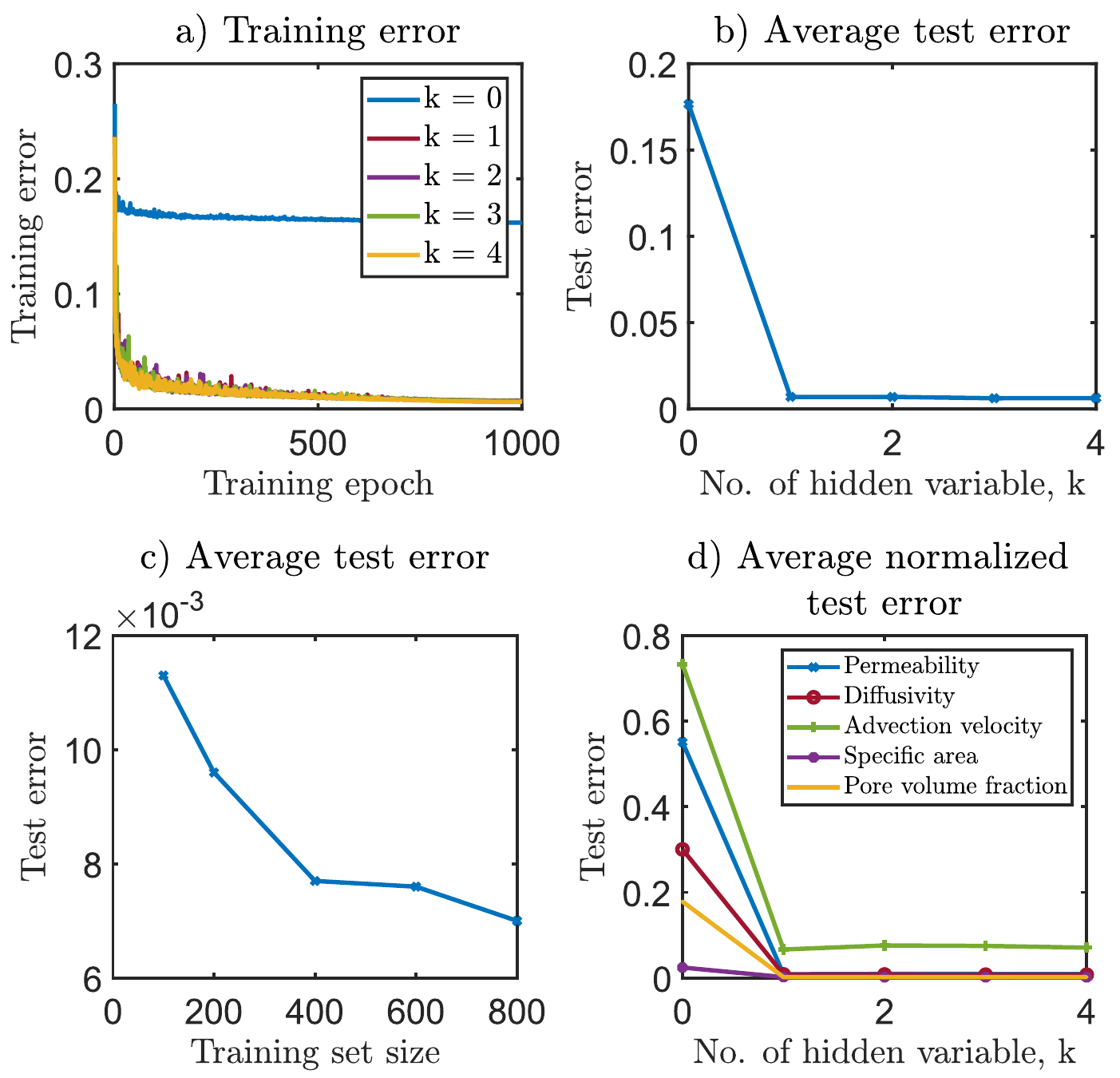}
    \caption{ Training and testing the RNO. (a) Training error vs.\ training epochs, (b) Average test loss of the trained RNOs vs.\ the number of hidden variables, (c) Average test loss of an RNO with one internal variable vs.\ training set size. (d) Average normalized test error vs.\ the number of internal variables.}
    \label{fig:loss}
\end{figure}

The results are shown in Figure \ref{fig:loss} in terms of the loss (\ref{eq:loss}).  Figure \ref{fig:loss}(a) shows how the training loss changes with the number of epochs with varying number of internal variables.   We see that an RNO with no internal variables (i.e., no history dependence) is unable to reduce the training error beyond a certain point.  However, RNOs with one or more internal variable can be trained to a high degree of accuracy.  Figure \ref{fig:loss}(b) shows the test loss (the same loss as (\ref{eq:loss}), but computed for the test data set) as a function of the number of internal variables for the trained RNO.  We see that the trained RNO with no internal variable provides a very poor approximation, but RNOs with one or more internal variable provides excellent approximation.  We repeat the training for various sizes of training data and the average test loss is shown in Figure \ref{fig:loss}(c) for the case of a single internal variable.  We see that the size (800) of our training data set is adequate, and the average test error is small in each component.

Figure \ref{fig:loss}(d) shows the average normalized test error of the various physical outputs.  For each physical output, we define the normalized test error as 
\begin{equation}
\text{average normalized test error in $\bfD^*$} = \left( \frac{1}{D_\text{test}} \sum_{d=1}^{D_\text{test}} 
\frac{\int^T_0\left|  (\bfD^*)^{\text{truth}}_d -  (\bfD^*)^{\text{approx}}_d \right|^2 d t}{\int^T_0\left|  (\bfD^*)^{\text{truth}}_d \right|^2 d t}\right)^{1/2}
\end{equation}
and so forth.  We observe that the average normalized test error about 1\% for each of the physical quantities except for the effective advection velocity $\bar{\bfv}$ where the error is about 6\% for a trained RNO with one or more internal variables.  The test includes cases where the effective advection velocity is zero up to machine precision and these lead to large apparent errors.

\begin{figure}[t]
\centering
    \includegraphics[width=0.7\textwidth]{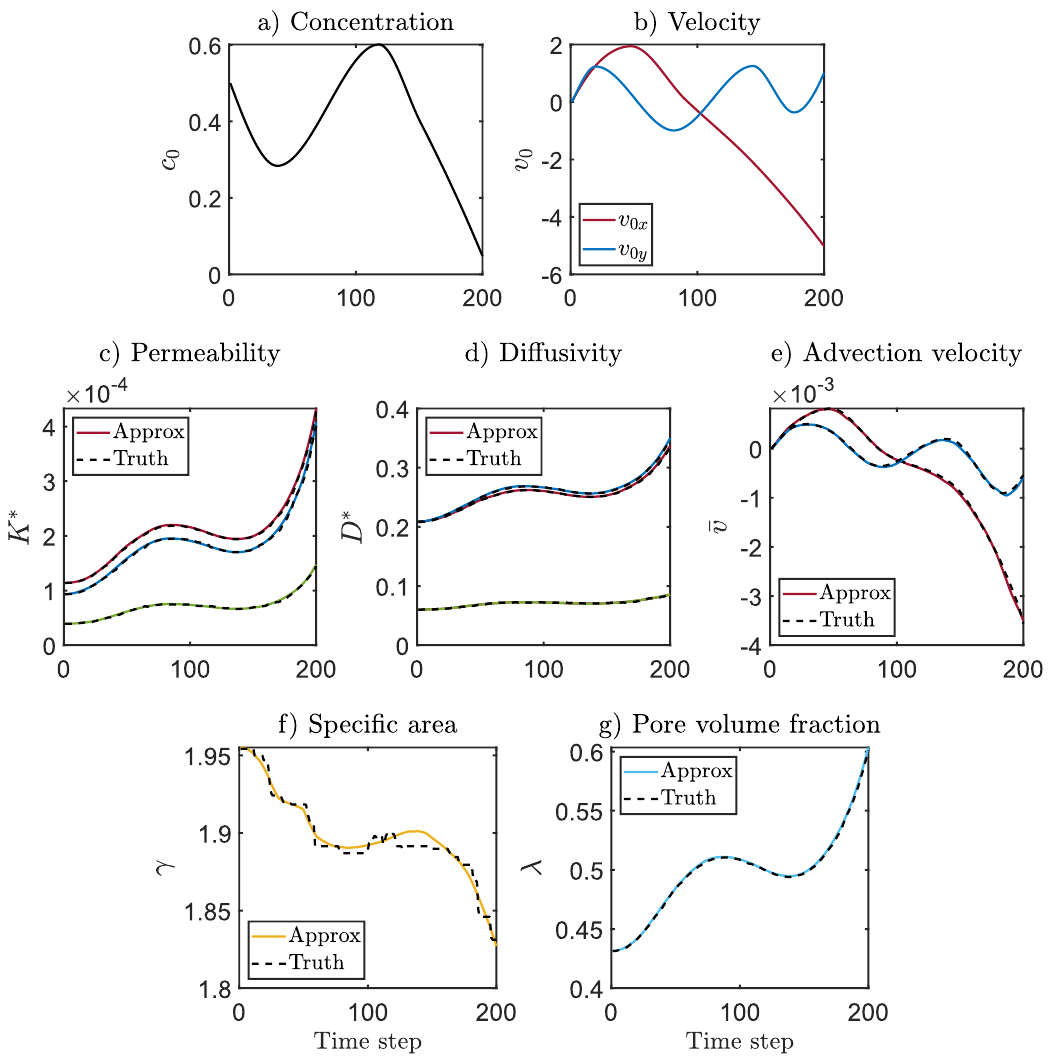}
    \caption{Input test trajectories, (a) concentration, (b) velocity trajectories,
    comparison of estimated values of the RNO with one hidden variable with ground truth, (c) permeability, (d) diffusivity, (e) advection velocity, (f) specific area and (g) pore volume fraction. }
    \label{fig:traj1}
\end{figure}

Figure \ref{fig:traj1} elaborates on the results by focussing on a typical trajectory chosen arbitrarily from the test data set.  Figures \ref{fig:traj1}(a,b) show the input, while Figures \ref{fig:traj1}(c-g) compare the ground truth and RNO predictions for the outputs: permeability ($\bfK^*$), diffusivity ($\bfD^*$), advection velocity ($\bar{\bfv}$), specific area ($\gamma$), and pore volume fraction ($\lambda$).  We see excellent agreement.

\begin{figure}
\centering
    \includegraphics[width=0.8\textwidth]{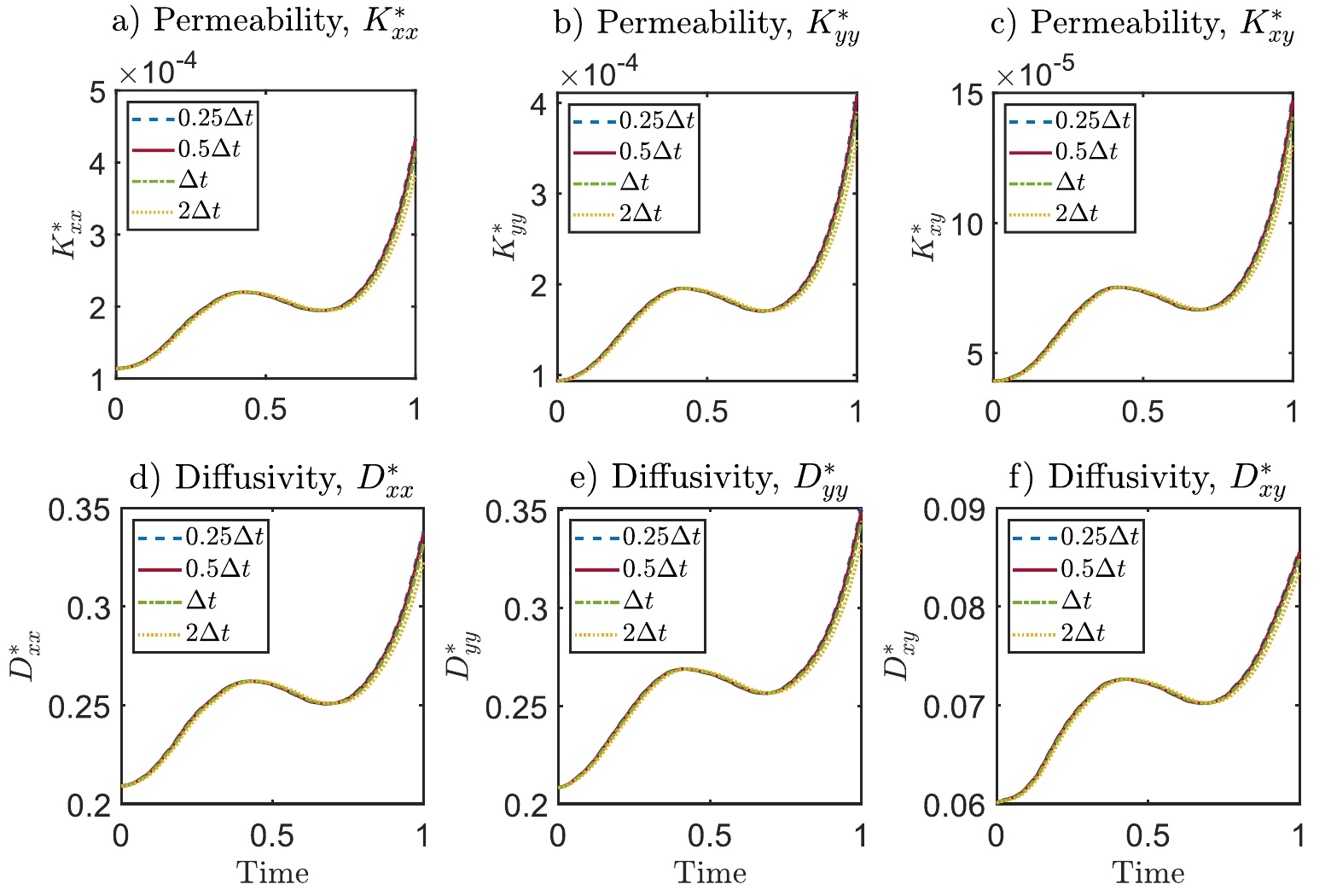}
    \caption{Comparison of estimated values of the RNO, permeability components (a)-(c), and diffusivity components (d)-(f) considering various time steps.}
    \label{fig:timeindependence}
\end{figure}

Finally, we examine the time discretization independence of the trained RNO. The RNO is initially trained with $\Delta t$.  Figure \ref{fig:timeindependence} shows the predictions of the RNO for permeability and diffusivity components evaluated with different time steps --  $0.25\Delta t$, $0.5\Delta t$, $\Delta t$, and $2\Delta t$ -- all for the same input trajectory.   We see that the results are independent of the time step. 
\vspace{0.1in}

In summary, we conclude that an RNO with one internal variable is able to provide an excellent approximation to the solution operator of the core scale model.

\section{Multiscale simulation} \label{sec:Multiscale simulation}

We now consider a geological scale simulation, but one that constantly updates the properties, from core scale calculations according to the framework described in Section \ref{sec:2scalemodel}.  However, instead of solving the core scale problem at each point at each instant, we use the trained neural approximation described in Section \ref{sec:RNO} as a surrogate for the core scale problem at each point at each instant.   Of particular interest is to understand how the formation and its properties, as well as the flow, change over long periods of time, and how such changes are magnified by heterogeneities in the formation.
We implement the geological scale problem with the Python finite elements library, FEniCS, and an unstructured mesh with triangular elements.  

\subsection{Reactive flow with uniform initial properties} \label{sec:uniform}
\begin{figure}
\centering
    \includegraphics[width=0.6\textwidth]{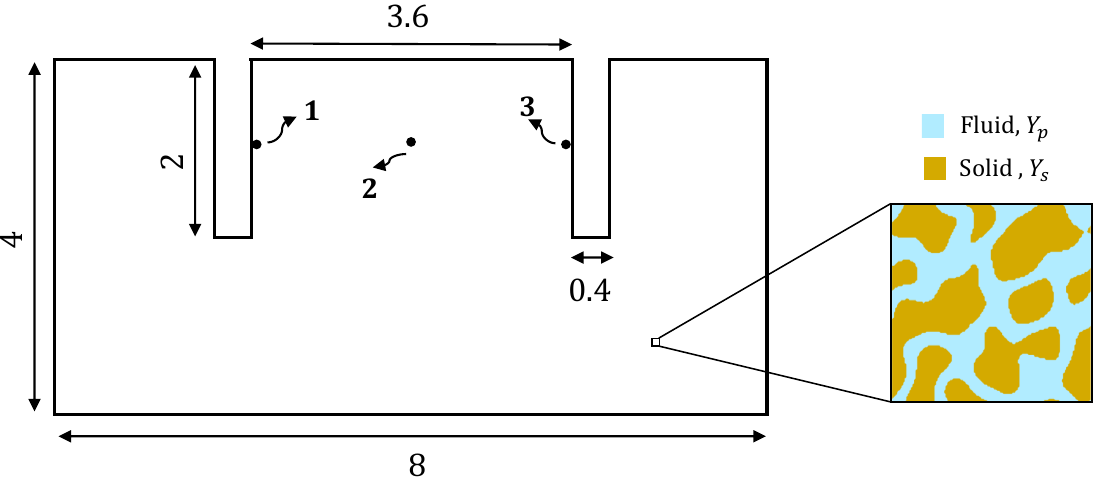}
    \caption{Geometry of geological domain.}
    \label{fig:Sec5_Fig1}
\end{figure}
We first consider an example with uniform initial material properties. The geometry of the geological formation and microstructure are shown in Figure \ref{fig:Sec5_Fig1}: a solution with high concentration of solute is injected at high pressure through the well on the left and removed from a well at the right.  The properties and boundary conditions are as follows.
\begin{align*}
    &\text{Pe} = 1000, \ \ \text{Da} = 0.001, \ \ c^* = 0.5,\\
    &\text{left well:} ~ p_0 = 10^5, ~~ c_0 = 0.6, \\
    &\text{right well:} ~ p_0 = 10^4, ~~ c_0 = 0.4, \\
    &\text{rest of the boundary: zero flux}.
\end{align*}

\begin{figure}
\centering
    \includegraphics[width=0.9\textwidth]{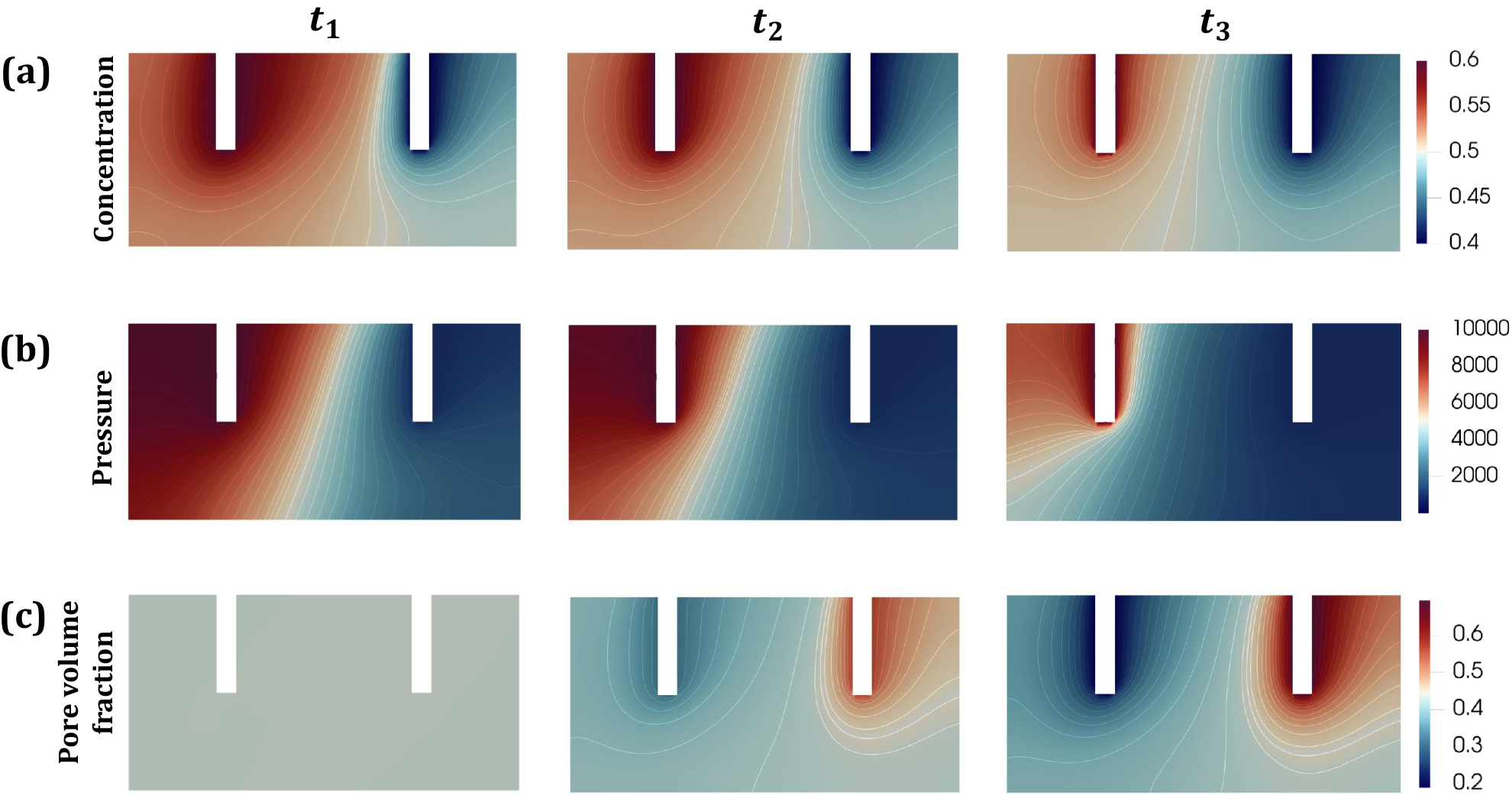}
    \caption{Variation of concentration (a)-(c), and pressure (d)-(f), and pore volume fraction (g)-(i) profiles in the geological formation due to chemical reactions. }
    \label{fig:Sec5_Fig2}
\end{figure}

Figures \ref{fig:Sec5_Fig2} shows three snapshots of the pressure and concentration profile at the geological scale, at $t_1$, $t_2$, and $t_3$  associated with time steps $1$, $110$, and $220$, respectively. The pressure and concentration change gradually away from the wells in the early times. However, at later times, the pressure and concentration change rapidly in the vicinity of the wells.  Recall that the equations (\ref{eq:reservoir}) at the geological scale are steady state equations, i.e., do not involve time explicitly.  Therefore, the evolution with time is related to the change in properties related to the reconstruction of the porous medium at the core scale. 

The concentration prescribed at the inlet exceeds the equilibrium concentration, resulting in the gradual precipitation of solutes on the solid/fluid interface in this region. Conversely, the concentration prescribed at the outlet lies below the equilibrium concentration, leading to a dissolution of solid structure close to the outlet. As illustrated in Figure \ref{fig:Sec5_Fig2}, precipitation in the inlet increases pore volume fraction that decreases permeability and diffusivity values while increasing the pressure gradient over time. This reduction in permeability and diffusivity reduces flow and chemical transport in the medium, eventually leading to clogging. 

\begin{figure}
\centering
    \includegraphics[height=8in]{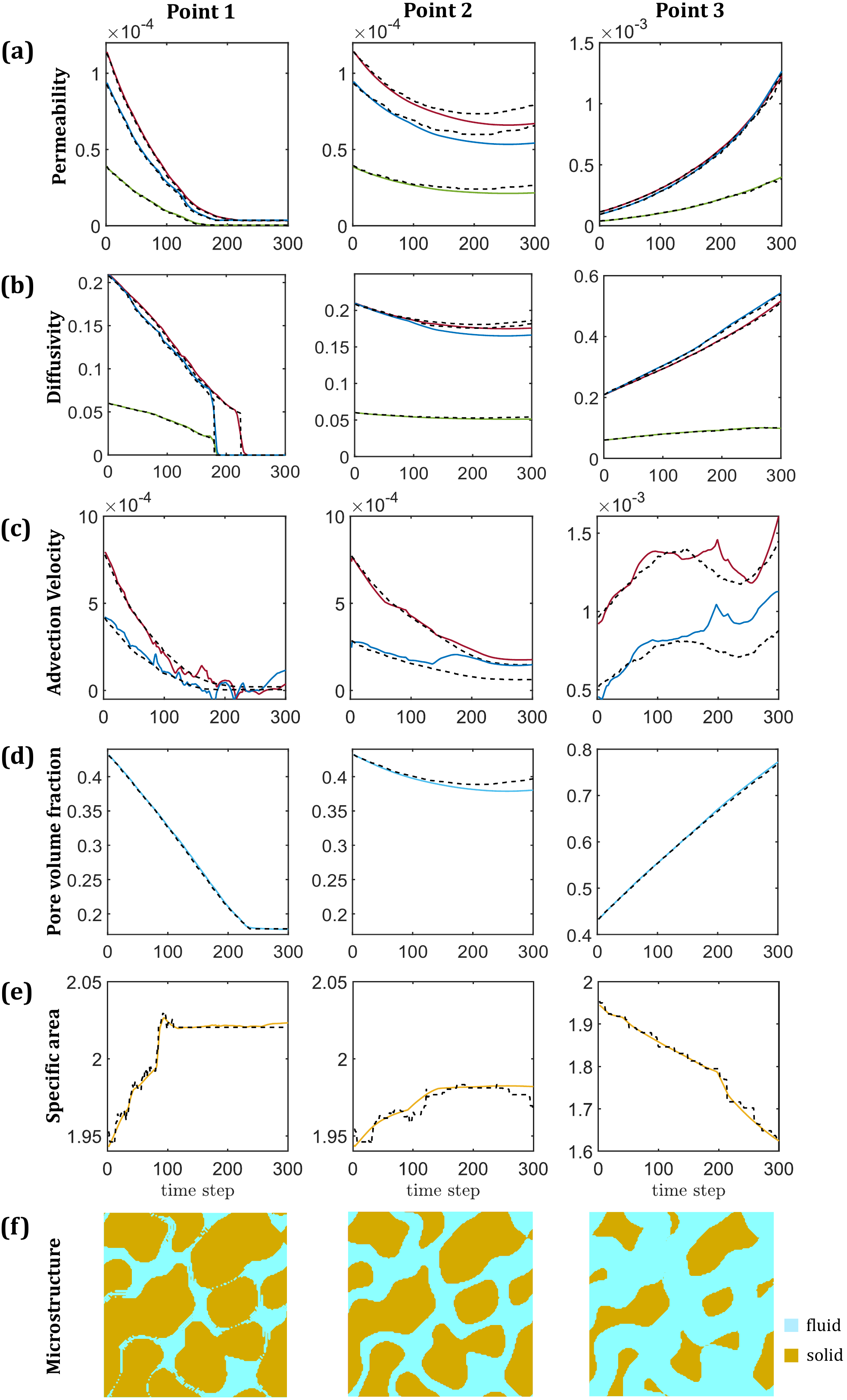}
    \caption{Comparison of the (a) permeability, (b) diffusivity, (c) advection velocity, (d) pore volume fraction, and (e) specific area components obtained from the RNO with ground truth, and (f) microstructure geometry at $230$th time step, at three random points.}
    \label{fig:comparison}
\end{figure}
We explore this further at the three points marked in Figure \ref{fig:Sec5_Fig1}: Points 1 and 3 are in the vicinity of the left well (inlet), and right well (outlet), respectively; Point 2 is located between the inlet and outlet. The changes in properties and pore structure with time are shown in Figure \ref{fig:comparison}  as the solid curves.  The left column of the figure shows the changes at Point 1 that is close to the inlet.  The deposition on the surface of solid leads to a  decrease in permeability, diffusivity, advection velocity and pore volume fraction.  The decrease is steady initially, but at some time the permeability, diffusivity and advection velocity effectively goes to zero: this time coincides with when the pores are blocked by deposition.  The porosity does not go to zero, but there is no change in its value beyond this time.  The specific surface area initially increases, but saturates earlier than full blocking.  The deposition leads to an increase of surface area, but this is eventually balanced by a decrease as the pores are blocked and opposite sides of a pore begin to touch.

We see the opposite trends on the right column corresponding to Point 3 that is close to the outlet.  We see that the dissolution leads to an increase of permeability, diffusivity, advection velocity and  pore volume fraction, and decrease of specific surface area.  The middle column shows the results for point 2 located between the inlet and outlet.  This only experiences slight changes in microstructure and effective properties, including permeability, diffusivity, pore volume fraction, and specific area. The reduction in the advection velocity at point 2 primarily results from the decrease in the flow reaching this point due to the clogging of the inlet.

Figure \ref{fig:comparison} provides an {\it a posteriori} evaluation of the error in using the RNO surrogate.  We take the velocity and concentration experienced by the 3 points during the macroscale calculations.  We compare the output of our RNO surrogate (solid lines) and the results of core scale calculations (dashed lines) for these histories.  We observe excellent overall agreement.  This tells us that our RNO surrogate performs well not only for the histories in the distribution used to train the model, but also real histories experienced in actual calculations.

\begin{figure}[t]
\centering
    \includegraphics[width=0.85\textwidth]{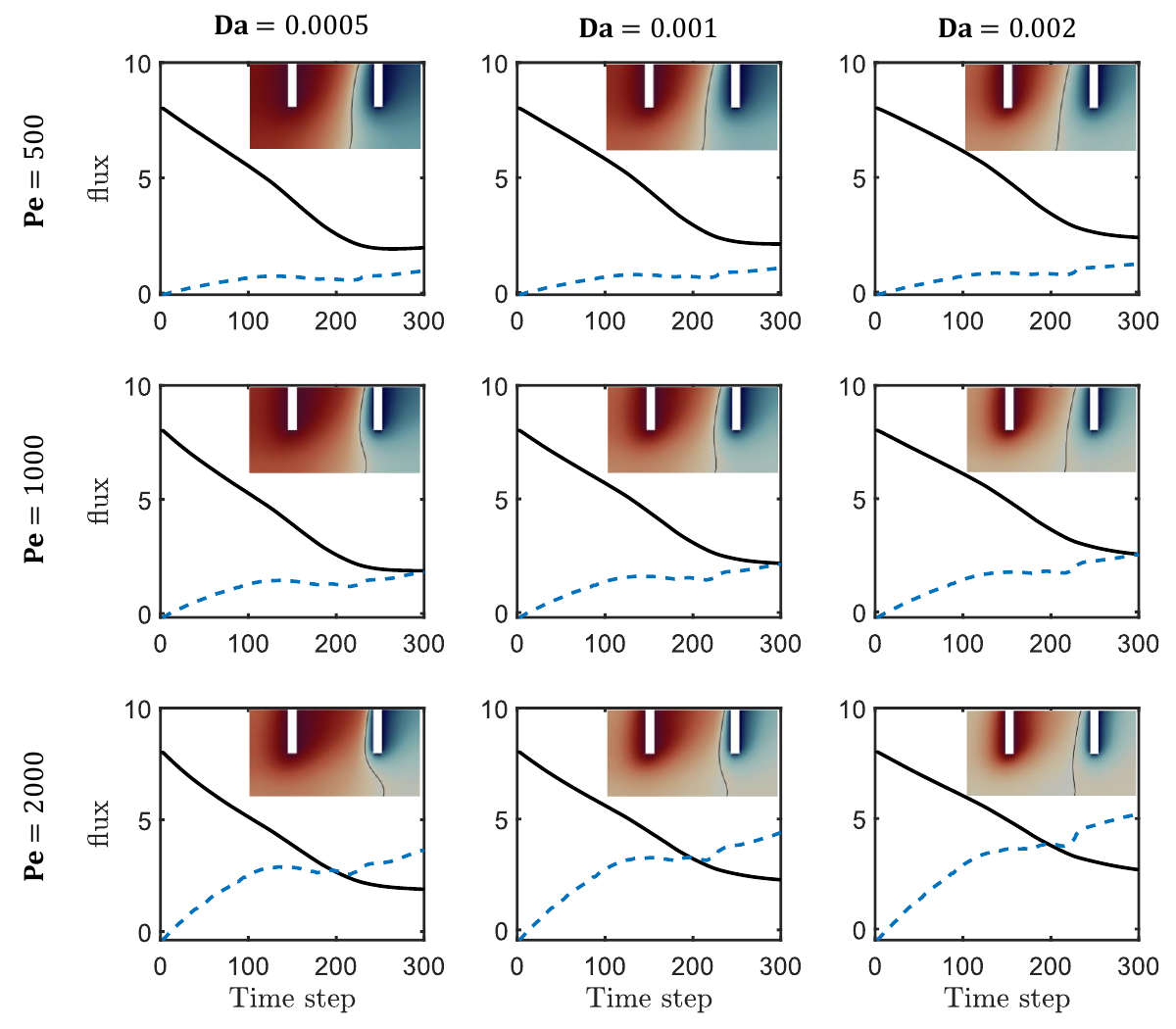}
    \caption{Variation of flow flux ($q$, solid line) in inlet and outlet, and chemical flux difference ($\Delta J$, dashed line) between inlet and outlet, considering various values of Peclet and Damkohler numbers. }
    \label{fig:flux}
\end{figure}

We now turn to the overall flow through the geometric formulation and the amount of material deposited in the formation as a function of time.  This is shown in the center of Figure \ref{fig:flux} for the parameters chosen above.  We see that the flow decreases and the amount of material deposited increases steadily till they saturate.  This is consistent with the observations above.  

Figure \ref{fig:flux} also shows a parameter study for various P\'eclet and Damk\"ohler numbers that reveal an interesting interplay between the pore and geological scale phenomena.  Recall that the core scale problem does not depend on the exact values of these non-dimensional quantities (as long as they satisfy the scaling).  It follows that the RNO is independent of them, and we do not have to retrain it for each pair of numbers.  Instead, these non-dimensional quantities only appear in the effective geological scale equations.  So, we repeat the geological scale calculation as described above for three distinct values (differing by a factor 2) of both these non-dimensional quantities.   We see that the flow decreases and the amount of material deposited increases steadily till they saturate in each case.   The change in flow and amount of material deposited increases significantly with increasing P\'eclet number, but only slightly with increasing Damk\"ohler number.  

Recall that we have very high P\'eclet number and a small Damk\"ohler number.  This means that we are in the reaction controlled regime.  Therefore, one could expect that the overall deposition rate would be more sensitive to the Damk\"ohler number, and less sensitive to the P\'eclet number.  However, we see the opposite trend in Figure \ref{fig:flux}.  This is because the morphological changes in the formation lead to clogging, and therefore the flow is transport limited at the macroscopic scale.  The inset in each graph in Figure \ref{fig:flux} shows a snapshot of the concentration at a fixed time: we see that the P\'eclet number induces a greater effect than the Damk\"ohler number.  In other words, even when the core scale is reaction limited, the formation can become transport limited.

\begin{table}
    \centering
    \caption{Comparison of computation cost (wall-clock time in seconds)}\label{tab:1}
    \begin{tabular}{ c c c } 
    \hline
    Method & Calculation & Cost \\
    \hline
    Classical constitutive relations & Geological scale calculation & 800 (CPU)) \\
    \hline
    Proposed method & Geological scale calculation & $900$ (CPU) \\ 
    & Off-line data (sequential) &  $3.6\times 10^6$ (CPU) \\
    & Off-line data (parallel) &  $4\times 10^4$ (GPU) \\
    & Off-line training & $5\times 10^3$ (GPU) \\
    \hline
    Concurrent multiscale (estimated) & Geological scale calculation & $1.64\times 10^7$ (CPU)) \\
    \hline
    \end{tabular}
\end{table}

Having established the accuracy and its efficacy in parameter study, we study the computational cost of the proposed approach and compare it with both the cost of the classical empirical constitutive model and the cost of concurrent multiscale models.  The results are shown in Table \ref{tab:1} for the base simulation above.  These calculations were conducted on a single core of an Intel Skylake CPU with a clock speed of 2.1 GHz and an NVIDIA P100 GPU with 3584 cores and 1.3 MHz clock speed. The classical empirical model is described in \cite{baqer2022review}. We find that the computational cost of solving the macroscopic problem using the trained RNO is comparable to the cost of classical constitutive relations.  These are significantly (by factor $10^5$) smaller than the estimated cost of the concurrent multiscale approach (we estimate this by using the cost of the unit cell problem and multiplying it by the product of the time steps and spatial grid).  The proposed approach has a one-time off-line cost of generating the data and fitting the RNO.  This is also smaller (by a factor $~ 10^2$) than solving a single concurrent multiscale calculation.  This off-line cost can further be reduced by parallelization.  
\vspace{0.1in}

In summary, we conclude that our approach is able to provide the accuracy of a concurrent multiscale model at the computational cost comparable to that of a classical constitutive relation.

\subsection{Reactive flow with non-uniform initial properties}
\begin{figure}
\centering
    \includegraphics[width=0.45\textwidth]{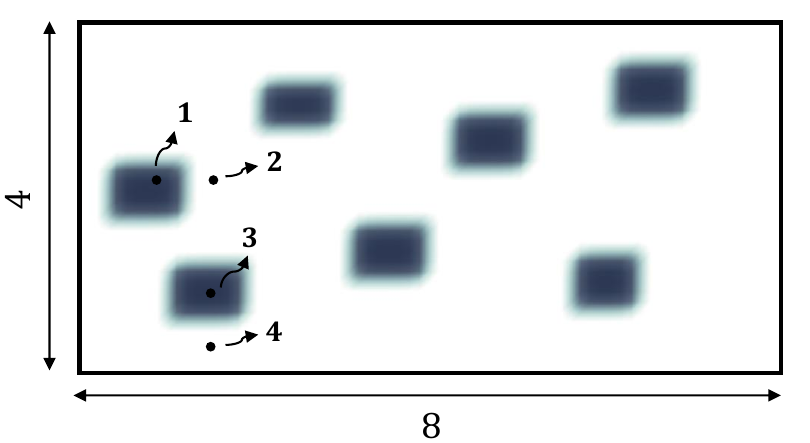}
    \caption{Domain geometry and location of blocks.}
    \label{fig:Sec4_Fig4}
\end{figure}

We now turn to examples with geological formation characterized by initially non-uniform properties as shown in Figure \ref{fig:Sec4_Fig4}.  This heterogeneity is defined by the initial value for the internal variable in the geological scale problem.

\subsubsection{High initial permeability and diffusivity inclusions} \label{sec:high}

\begin{figure}
\centering
    \includegraphics[width=0.9\textwidth]{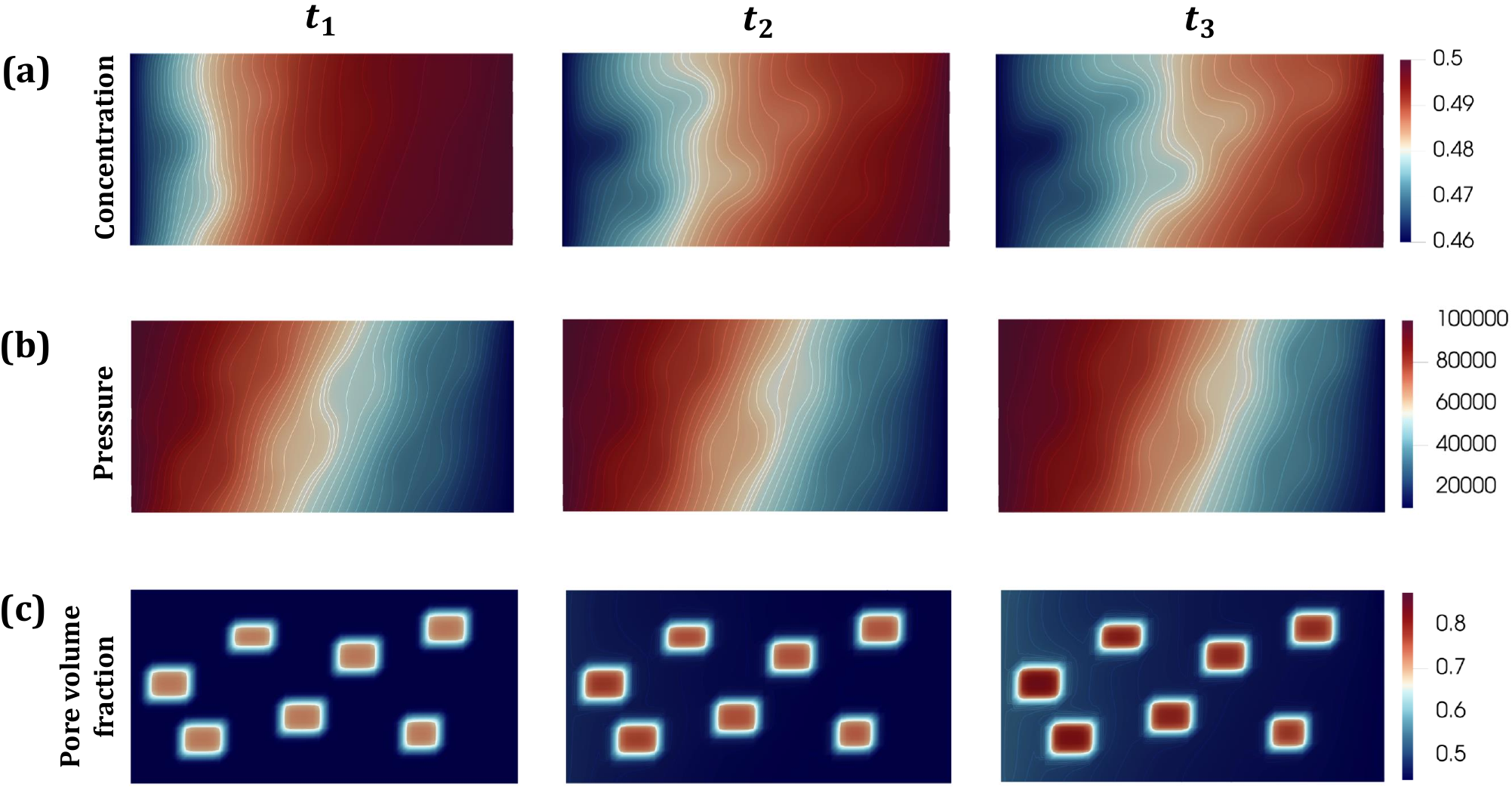}
    \caption{Variation of concentration (a), pressure (b), and pore volume fraction (c) fields, considering non-uniform initial material properties.}
    \label{fig:Sec5_Fig5}
\end{figure}

We consider a $4\times 8 ~ \text{m}$ domain with embedded blocks of higher initial permeability and diffusivity in the geometry illustrated in Figure \ref{fig:Sec4_Fig4}. A horizontal flow is introduced from the left boundary and withdrawn from the right.  The parameters are the same as in Section \ref{sec:uniform}, and boundary conditions are
\begin{align*}
    &\text{left boundary:} ~ p_0 = 10^5, ~~ c_0 = 0.4, \\
    &\text{right boundary:} ~ p_0 = 10^4, ~~ c_0 = 0.5, \\
    &\text{rest of the boundary: zero flux}.
\end{align*}
Figures \ref{fig:Sec5_Fig5} (a), and (b) show the temporal evolution of the concentration and pressure fields over time due to changes in the morphology.  The left (inlet) concentration is lower than the equilibrium concentration, and this leads to dissolution at that end.  The right (outlet) concentration is at the equilibrium concentration, and so we do not see significant morphological changes there.  Note that the flow and concentration are not uniform across the (vertical) cross-section even though the boundary conditions are.  This is a result of the inclusions.  The flow preferentially enters the regions with high permeability and diffusivity, and seeks to connect these regions together.  The greater flow leads to greater chemical reaction and further increases permeability; this is clear from the evolution of the pore volume fraction shown in Figure \ref{fig:Sec5_Fig5} (c)).  This further aids the channeling from one inclusion to another.

\begin{figure}
\centering
    \includegraphics[width=0.8\textwidth]{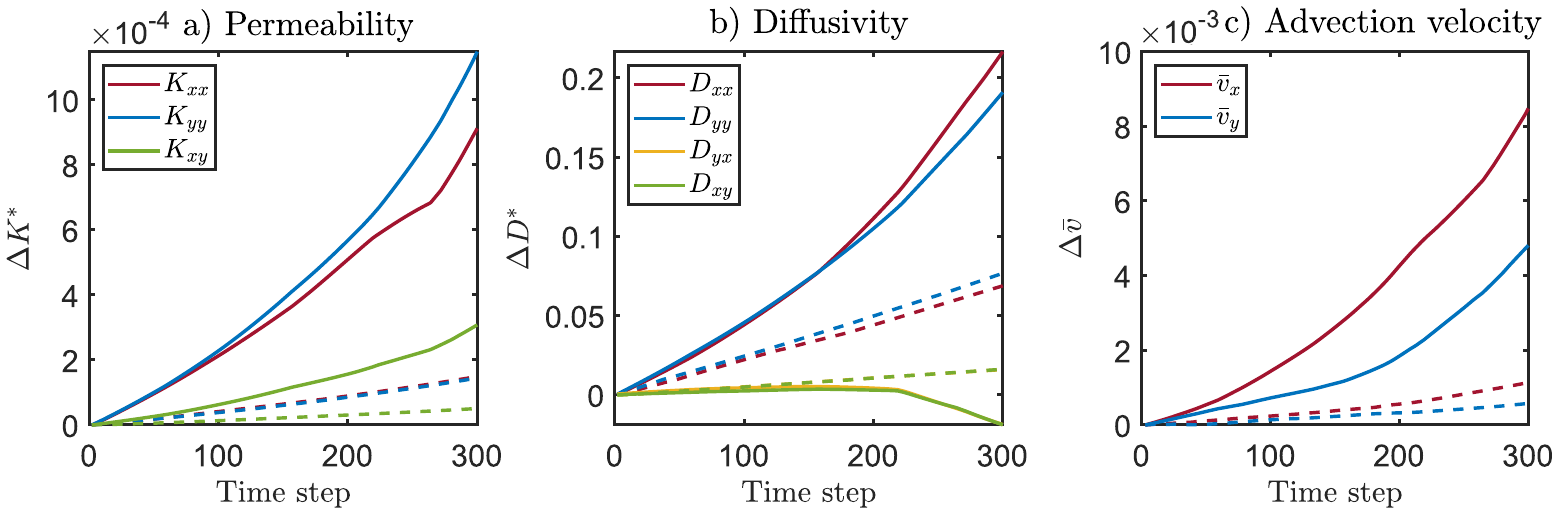}
    \caption{Change in effective properties, (a) permeability, (b) diffusivity, and (c) advection velocity, at points 1 (solid lines), and point 2 (dashed lines).}
    \label{fig:Sec4_Fig6}
\end{figure}

The differential change in permeability, diffusivity, and advection velocity is emphasized in Figure \ref{fig:Sec4_Fig6} that shows the evolution of these quantities at two points  (marked 1 and 2 in Figure \ref{fig:Sec4_Fig4}).  The solid lines correspond to Point 1 located within the high permeability block, while the dashed lines represent the change of effective properties at Point 2 located outside the high permeability block. As expected, the rate of change in effective properties at Point 1 is significantly greater compared to Point 2. 

\subsubsection{Low initial permeability and diffusivity inclusions}\label{sec:low}

\begin{figure}
\centering
    \includegraphics[width=0.9\textwidth]{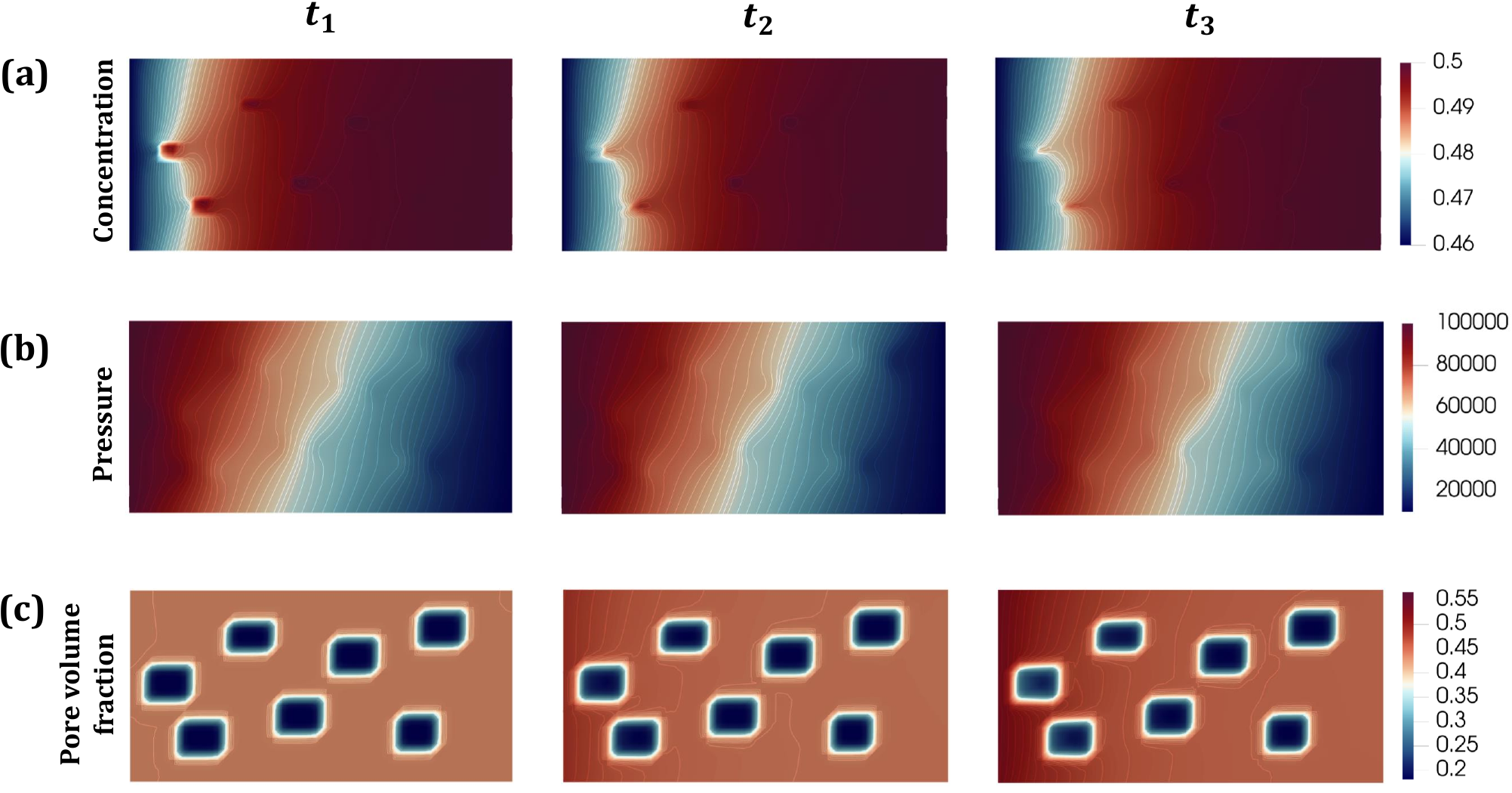}
    \caption{Variation of concentration (a), pressure (b), and pore volume fraction (c) fields, considering non-uniform initial material properties.}
    \label{fig:Sec5_Fig6}
\end{figure}

\begin{figure}
\centering
    \includegraphics[width=0.8\textwidth]{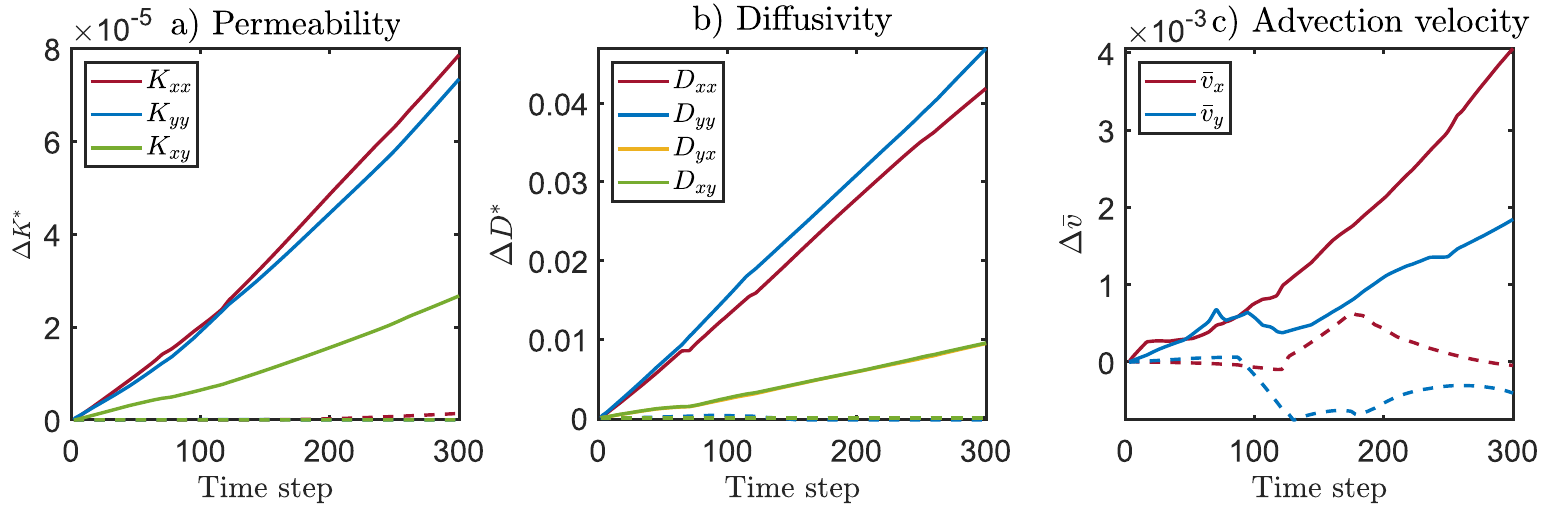}
    \caption{Change in effective properties, (a) permeability, (b) diffusivity, and (c) advection velocity, at points 3 (dashed lines), and point 4 (solid lines).}
    \label{fig:Sec4_Fig8}
\end{figure}

We now consider the complementary situation with blocks of lower initial permeability and diffusivity.   The parameters are the same as in Section \ref{sec:uniform}, and boundary conditions are the same as Section \ref{sec:high}.    Figure \ref{fig:Sec5_Fig6} (a)-(b) shows that the flow moves around the block following regions with higher initial permeability and diffusivity.  Consequently, there is more reaction and dissolution in the matrix compared to the blocks, Figure \ref{fig:Sec5_Fig6}(c).  This further channels the flow through between the blocks.  Figure \ref{fig:Sec4_Fig8} shows the change in material properties at Points 3 and 4 (as shown in figure \ref{fig:Sec4_Fig4}). We observe a substantial change in effective properties at Point 4 outside the blocks and relatively little change at Point 3 inside.  

\section{Conclusion} \label{sec:conc}

The transport of water and solutes through permeable geological formations is a complex multiscale phenomenon.  In this work, we have proposed a framework that harnesses the power of continuum modeling and machine learning to address this complexity with accuracy at reasonable computational cost.  
We have demonstrated the framework on flow through permeable geological formations with advection-dominated transport and volume reactions.

We begin with an asymptotic analysis of the governing equations that enable us to split the problem into two interconnected problems: one at the scale of the geological formation and another at the scale of a core.  A key idea in this analysis is the invocation of a drifting coordinate system to capture the core.  We then introduce a recurrent neural operator (RNO) to approximate the solution operator of the core-scale problem.  This consists of two feed-forward neural networks and internal variables.  The neural networks are trained and the internal variables are identified from data that is generated by repeatedly solving the core-scale problem.  The key features of this neural architecture are that it is consistent with common formulations of continuum physics, that it is relatively simple and that it is independent of the time-discretization.  We demonstrate that it is able to accurately capture the behavior of the core scale over long periods of time including the morphological changes in the pores and the resulting change in effective permeability, diffusivity, advective velocity, pore volume fraction and specific area.  Finally, we solve the problem of transport and morphological changes at the scale of the geological formation by using the trained RNO as a surrogate for the small scale problem.  We thus obtain the accuracy of a concurrent multiscale simulations at a cost that is comparable to classical constitutive relations.  We demonstrate the ability of this approach to learn subtle features of the interaction between the scales including the change from reaction limited to transport limited regime due to clogging and the positive feedback of channeling in heterogeneous situations.

We now emphasize a few notable aspects of the proposed approach.  First, our RNO neural approximation is able to capture morphological and property changes over long periods of time.  It is formulated in a time-continuous manner and discretized as necessary for training and use.  It follows that the approximation is independent of the discretization.  This is important because it is common to use different time discretization for the core scale problem used to generate data and the application at the geological scale.  It also enables the use of data generated at different discretizations.  Second, the two scale formulation makes the core scale problem to be independent of the physical characteristics of the flow and reaction rate, specifically the P\'eclet and the Damk\"ohler numbers (as long as we are in the advection dominated regime with slow reactions).  This means that we need to generate data and train the RNO only once, and can use the method for different situations as these quantities change.  Third, we not only demonstrate accuracy over the distribution used to train the RNO, but also {\it a posteriori} over the actual histories encountered in the geological scale calculations.  Finally, the approach is highly transferable.  In this work, we used examples in two dimensions to demonstrate the framework with modest computational cost.  However, the approach holds in three dimensions. Similarly, one can incorporate other phenomena including, for example, multi-phase flows, more complex chemistry, poroelasticity and phase change (melting), as long as we can use scale separation.  Importantly, one can extend this to more than two scales as long as they interact pairwise.  

We close with a discussion of a few limitations and avenues for future work.  First, the approach requires us to train the RNO with a starting pore morphology.  This is not an issue as long as the core scale is chosen to be large enough to be statistically representative of the underlying medium.  However, this adds to the computational cost.  So, one may consider training the RNO on an ensemble of cores.   Second, while we demonstrate {\it a posteriori} accuracy, it would be useful to have a systematic approach to quantifying the overall uncertainties. Such a quantification may also enable the use of active learning where we initially train the RNO over synthetic samples as we do in this work, but then progressively add more data based on examples we encounter.  Third, in this work, we only use geological scale information obtained by averaging the results of core scale simulations to train the RNO.  However, we have access to core scale information.  It is possible that this data has insights that may lead to a more robust efficient training procedure.  Fourth, we have exclusively used data generated numerically to train the RNO.  it would be interesting to use a combination of experimental and computational data.  Finally, it remains to use the framework established in this work to actual geological problems.  All of this remains a topic of current and future research.

\section*{Data availability}
The data is available at CaltechData through \url{https://doi.org/10.22002/yd0c5-q5s36}  

\section*{Acknowledgments}
We are delighted to acknowledge numerous helpful discussions with Professor Xiaojing (Ruby) Fu.  We gratefully acknowledge the financial support of the Resnick Sustainability Institute at the California Institute of Technology.  KB also acknowledges the support of the Army Research Office through grant number W911NF-22-1-0269.  The simulations reported here were conducted on the Resnick High Performance Computing Cluster at the California Institute of Technology.

\bibliographystyle{abbrv}
\bibliography{poromech-refs}

\begin{thebibliography}{10}

\bibitem{allaire2010two}
G.~Allaire, R.~Brizzi, A.~Mikeli{\'c}, and A.~Piatnitski.
\newblock Two-scale expansion with drift approach to the taylor dispersion for reactive transport through porous media.
\newblock {\em Chemical Engineering Science}, 65:2292--2300, 2010.

\bibitem{allaire2007homogenization}
G.~Allaire and A.-L. Raphael.
\newblock Homogenization of a convection--diffusion model with reaction in a porous medium.
\newblock {\em Comptes Rendus Mathematique}, 344:523--528, 2007.

\bibitem{appelo2004geochemistry}
C.~A.~J. Appelo and D.~Postma.
\newblock {\em Geochemistry, groundwater and pollution}.
\newblock CRC Press, 2004.

\bibitem{auriault1995taylor}
J.-L. Auriault and P.~Adler.
\newblock Taylor dispersion in porous media: analysis by multiple scale expansions.
\newblock {\em Advances in Water Resources}, 18:217--226, 1995.

\bibitem{baqer2022review}
Y.~Baqer and X.~Chen.
\newblock A review on reactive transport model and porosity evolution in the porous media.
\newblock {\em Environmental Science and Pollution Research}, 29:47873--47901, 2022.

\bibitem{bear2012introduction}
J.~Bear and Y.~Bachmat.
\newblock {\em Introduction to modeling of transport phenomena in porous media}.
\newblock Springer Science \& Business Media, 2012.

\bibitem{bhattacharya_2023}
K.~Bhattacharya, B.~Liu, A.~Stuart, and M.~Trautner.
\newblock Learning markovian homogenized models in viscoelasticity.
\newblock {\em Multiscale Modeling \& Simulation}, 21:641--679, 2023.

\bibitem{carbonell1983dispersion}
R.~Carbonell and S.~Whitaker.
\newblock Dispersion in pulsed systems. {II}: {T}heoretical developments for passive dispersion in porous media.
\newblock {\em Chemical Engineering Science}, 38:1795--1802, 1983.

\bibitem{demirer2022improving}
E.~Demirer, E.~Coene, A.~Iraola, A.~Nardi, E.~Abarca, A.~Idiart, G.~de~Paola, and N.~Rodr{\'\i}guez-Morillas.
\newblock Improving the performance of reactive transport simulations using artificial neural networks.
\newblock {\em Transport in Porous Media}, 149:1--27, 2022.

\bibitem{donato2005averaging}
P.~Donato and A.~Piatnitsk\'i.
\newblock Averaging of nonstationary parabolic operators with large lower order terms.
\newblock In {\em Multi Scale Problems and Asymptotic Analysis}, pages 153--165. Gakkotosho, 2006.

\bibitem{ghavamian2019accelerating}
F.~Ghavamian and A.~Simone.
\newblock Accelerating multiscale finite element simulations of history-dependent materials using a recurrent neural network.
\newblock {\em Computer Methods in Applied Mechanics and Engineering}, 357:112594, 2019.

\bibitem{kang2003simulation}
Q.~Kang, D.~Zhang, and S.~Chen.
\newblock Simulation of dissolution and precipitation in porous media.
\newblock {\em Journal of Geophysical Research: Solid Earth}, 108:2505, 2003.

\bibitem{kang2002lattice}
Q.~Kang, D.~Zhang, S.~Chen, and X.~He.
\newblock Lattice {B}oltzmann simulation of chemical dissolution in porous media.
\newblock {\em Physical Review E}, 65:036318, 2002.

\bibitem{karimi_2023}
M.~Karimi and K.~Bhattacharya.
\newblock A fast-{F}ourier transform method for reactive flow in porous media.
\newblock {\em In preparation}, 2023.

\bibitem{kinga2015method}
D.~Kingma and J.~Ba.
\newblock {ADAM: A} method for stochastic optimization.
\newblock In {\em International conference on learning representations (ICLR)}, volume~5, page~6, 2015.

\bibitem{klambauer2017self}
G.~Klambauer, T.~Unterthiner, A.~Mayr, and S.~Hochreiter.
\newblock Self-normalizing neural networks.
\newblock {\em Advances in Neural Information Processing Systems}, 30, 2017.

\bibitem{levy1983fluid}
T.~L{\'e}vy.
\newblock Fluids in porous media: {D}arcy's law.
\newblock In {\em Homogenization Techniques for Composite Media}, pages 75--84. Springer, 1985.

\bibitem{li2008level}
X.~Li, H.~Huang, and P.~Meakin.
\newblock Level set simulation of coupled advection-diffusion and pore structure evolution due to mineral precipitation in porous media.
\newblock {\em Water Resources Research}, 44, 2008.

\bibitem{lichtner2018reactive}
P.~C. Lichtner, C.~I. Steefel, and E.~H. Oelkers.
\newblock {\em Reactive transport in porous media}, volume~34.
\newblock Walter de Gruyter GmbH \& Co KG, 2018.

\bibitem{liu2022learning}
B.~Liu, N.~Kovachki, Z.~Li, K.~Azizzadenesheli, A.~Anandkumar, A.~M. Stuart, and K.~Bhattacharya.
\newblock A learning-based multiscale method and its application to inelastic impact problems.
\newblock {\em Journal of the Mechanics and Physics of Solids}, 158:104668, 2022.

\bibitem{liu2023learning}
B.~Liu, E.~Ocegueda, M.~Trautner, A.~M. Stuart, and K.~Bhattacharya.
\newblock Learning macroscopic internal variables and history dependence from microscopic models.
\newblock {\em Journal of the Mechanics and Physics of Solids}, 178:105329, 2023.

\bibitem{liu2022machine}
M.~Liu, B.~Kwon, and P.~K. Kang.
\newblock Machine learning to predict effective reaction rates in 3d porous media from pore structural features.
\newblock {\em Scientific Reports}, 12:5486, 2022.

\bibitem{lu2021data}
H.~Lu, D.~Ermakova, H.~M. Wainwright, L.~Zheng, and D.~M. Tartakovsky.
\newblock Data-informed emulators for multi-physics simulations.
\newblock {\em Journal of Machine Learning for Modeling and Computing}, 2(2), 2021.

\bibitem{maruvsic2005homogenization}
E.~Maru{\v{s}}i{\'c}-Paloka and A.~L. Piatnitski.
\newblock Homogenization of a nonlinear convection-diffusion equation with rapidly oscillating coefficients and strong convection.
\newblock {\em Journal of the London Mathematical Society}, 72:391--409, 2005.

\bibitem{mauri1991dispersion}
R.~Mauri.
\newblock Dispersion, convection, and reaction in porous media.
\newblock {\em Physics of Fluids A: Fluid Dynamics}, 3:743--756, 1991.

\bibitem{medici_2021}
G.~Medici, L.~Smeraglia, A.~Torabi, and C.~Botter.
\newblock Review of modeling approaches to groundwater flow in deformed carbonate aquifers.
\newblock {\em Groundwater}, 59:334--351, 2021.

\bibitem{medsker2001recurrent}
L.~R. Medsker and L.~Jain.
\newblock Recurrent neural networks.
\newblock {\em Design and Applications}, 5(64-67):2, 2001.

\bibitem{mozaffar2019deep}
M.~Mozaffar, R.~Bostanabad, W.~Chen, K.~Ehmann, J.~Cao, and M.~Bessa.
\newblock Deep learning predicts path-dependent plasticity.
\newblock {\em Proceedings of the National Academy of Sciences}, 116:26414--26420, 2019.

\bibitem{oron1998flow}
A.~P. Oron and B.~Berkowitz.
\newblock Flow in rock fractures: The local cubic law assumption reexamined.
\newblock {\em Water Resources Research}, 34:2811--2825, 1998.

\bibitem{paine1983dispersion}
M.~Paine, R.~Carbonell, and S.~Whitaker.
\newblock Dispersion in pulsed systems—i: Heterogenous reaction and reversible adsorption in capillary tubes.
\newblock {\em Chemical Engineering Science}, 38:1781--1793, 1983.

\bibitem{quintard1993transport}
M.~Quintard and S.~Whitaker.
\newblock Transport in ordered and disordered porous media: volume-averaged equations, closure problems, and comparison with experiment.
\newblock {\em Chemical Engineering Science}, 48:2537--2564, 1993.

\bibitem{rubinstein1986dispersion}
J.~Rubinstein and R.~Mauri.
\newblock Dispersion and convection in periodic porous media.
\newblock {\em SIAM Journal on Applied Mathematics}, 46:1018--1023, 1986.

\bibitem{seigneur2019reactive}
N.~Seigneur, K.~U. Mayer, and C.~I. Steefel.
\newblock Reactive transport in evolving porous media.
\newblock {\em Reviews in Mineralogy and Geochemistry}, 85:197--238, 2019.

\bibitem{sprocati2021integrating}
R.~Sprocati and M.~Rolle.
\newblock Integrating process-based reactive transport modeling and machine learning for electrokinetic remediation of contaminated groundwater.
\newblock {\em Water Resources Research}, 57:e2021WR029959, 2021.

\bibitem{steefel1998multicomponent}
C.~I. Steefel and P.~C. Lichtner.
\newblock Multicomponent reactive transport in discrete fractures: I. controls on reaction front geometry.
\newblock {\em Journal of Hydrology}, 209:186--199, 1998.

\bibitem{sun2006dynamic}
S.~Sun and M.~F. Wheeler.
\newblock A dynamic, adaptive, locally conservative, and nonconforming solution strategy for transport phenomena in chemical engineering.
\newblock {\em Chemical Engineering Communications}, 193:1527--1545, 2006.

\bibitem{tartakovsky2016smoothed}
A.~M. Tartakovsky, N.~Trask, K.~Pan, B.~Jones, W.~Pan, and J.~R. Williams.
\newblock Smoothed particle hydrodynamics and its applications for multiphase flow and reactive transport in porous media.
\newblock {\em Computational Geosciences}, 20:807--834, 2016.

\bibitem{taylor1953dispersion}
G.~I. Taylor.
\newblock Dispersion of soluble matter in solvent flowing slowly through a tube.
\newblock {\em Proceedings of the Royal Society of London. Series A. Mathematical and Physical Sciences}, 219:186--203, 1953.

\bibitem{varloteaux2013reactive}
C.~Varloteaux, M.~T. Vu, S.~B{\'e}kri, and P.~M. Adler.
\newblock Reactive transport in porous media: pore-network model approach compared to pore-scale model.
\newblock {\em Physical Review E}, 87:023010, 2013.

\bibitem{wang2021physics}
K.~Wang, Y.~Chen, M.~Mehana, N.~Lubbers, K.~C. Bennett, Q.~Kang, H.~S. Viswanathan, and T.~C. Germann.
\newblock A physics-informed and hierarchically regularized data-driven model for predicting fluid flow through porous media.
\newblock {\em Journal of Computational Physics}, 443:110526, 2021.

\bibitem{wang2021upscaling}
Z.~Wang and I.~Battiato.
\newblock Upscaling reactive transport and clogging in shale microcracks by deep learning.
\newblock {\em Water Resources Research}, 57(4):e2020WR029125, 2021.

\bibitem{wang2022pore}
Z.~Wang, L.~Chen, H.~Wei, Z.~Dai, Q.~Kang, and W.-Q. Tao.
\newblock Pore-scale study of mineral dissolution in heterogeneous structures and deep learning prediction of permeability.
\newblock {\em Physics of Fluids}, 34, 2022.

\bibitem{witherspoon1980validity}
P.~A. Witherspoon, J.~S. Wang, K.~Iwai, and J.~E. Gale.
\newblock Validity of cubic law for fluid flow in a deformable rock fracture.
\newblock {\em Water Resources Research}, 16:1016--1024, 1980.

\bibitem{wu2020recurrent}
L.~Wu, N.~G. Kilingar, L.~Noels, et~al.
\newblock A recurrent neural network-accelerated multi-scale model for elasto-plastic heterogeneous materials subjected to random cyclic and non-proportional loading paths.
\newblock {\em Computer Methods in Applied Mechanics and Engineering}, 369:113234, 2020.

\end{thebibliography}

\newpage
\renewcommand\thefigure{S\arabic{figure} \Alph{section}}
\setcounter{figure}{0}
\renewcommand\thepage{S\arabic{page}}
\setcounter{page}{1}

\end{document}